\input ./style/arxiv-general.cfg
\documentclass[MSNbibl,number,citesort,seceqn,dvips]{arxbj}
\makeatletter
   \@ifpackageloaded{graphicx}{}{\usepackage{graphicx}}
\makeatother
\volume{22}
\issue{3}
\pubyear{2016}
\firstpage{1491}
\lastpage{1519}
\doi{10.3150/15-BEJ700}% Updated by VTEXPTS2LaTeX.exe, 29.04.2015 10:47
\docsubty{FLA}

\makeatletter
\def\mid{|}
\def\overset{\stackrel}
\def\cal{\mathcal}

\newcommand{\eqref}[1]{(\ref{#1})}
\newtheorem{thmm}{Theorem}[section]
\newtheorem{proposition}{Proposition}[section]
\newtheorem{lemma}{Lemma}[section]
\newtheorem{corollary}{Corollary}[section]
\newremark{rem}{Remark}[section]

\newcommand{\R}{\mathbb{R}}
\newcommand{\veps}{\varepsilon}
\newcommand{\cH}{{\cal H}}
\newcommand{\cY}{{\cal Y}}
\newcommand{\topr}{\stackrel{\mathrm{P}}{\longrightarrow}}
\newcommand{\RV}{\operatorname{RV}}
\makeatother

\begin{document}
\begin{frontmatter}

\title{Passage time and fluctuation calculations for subexponential L\'
{e}vy processes}
\runtitle{Fluctuations of subexponential L\'evy processes}

\begin{aug}
%%%% inicialai - be tarpu
% Corresponding author: Claudia Kluppelberg - cklu@tum.de% Updated by
%VTEXPTS2LaTeX.exe, 29.04.2015 10:47
\author[A]{\inits{R.}\fnms{Ron} \snm{Doney}\thanksref{A}\ead[label=e1]{rad@maths.manchester.ac.uk}},
\author[B]{\inits{C.}\fnms{Claudia} \snm{Kl\"uppelberg}\corref{}\thanksref{B}\ead[label=e2]{cklu@tum.de}}
\and
\author[C]{\inits{R.}\fnms{Ross} \snm{Maller}\thanksref{C}\ead[label=e3]{Ross.Maller@anu.edu.au}}
%\author[A]{\inits{}\fnms{}~\snm{}\corref{}\thanksref{A}
%\ead[label=e1]{}}%,
%\author[]{\inits{}\fnms{}~\snm{}\thanksref{}\ead[label=]{}}
% \and
%\author[]{\inits{}\fnms{}~\snm{}\thanksref{}\ead[label=]{}}
%%\runauthor{} %% auto
%\dedicated{}
\address[A]{School of Mathematics, Manchester University, Oxford Rd,
Manchester M13 9PL, United Kingdom.\\
\printead{e1}}

\address[B]{Center for Mathematical Sciences, Technische Universit\"at
M\"unchen,
Boltzmannstrasse 3, 85748 Garching, Germany.
\printead{e2}}

\address[C]{Research School of Finance, Actuarial Studies and Applied Statistics,
Australian National University, Canberra, ACT,  0200, Australia.
\printead{e3}}
%\address[A]{. \printead{e1}}
%\address[]{. \printead{}}
\end{aug}

% HISTORY:
%
\received{\smonth{4} \syear{2014}}% Updated by VTEXPTS2LaTeX.exe,
%29.04.2015 10:47
%
\revised{\smonth{10} \syear{2014}}% Updated by VTEXPTS2LaTeX.exe,
%29.04.2015 10:47

% ABSTRACT
%
\begin{abstract}
We consider the passage time problem for L\'evy processes, emphasising heavy
tailed cases. %%, with a view to applications in insurance risk.
Results are obtained under quite mild assumptions, namely, drift to
$-\infty$
a.s. of the process, possibly at a linear rate (the finite mean case), but
possibly much faster (the infinite mean case), together with
subexponential growth
on the positive side.
Local and functional versions of limit distributions are derived
for the passage time itself, as well as for the position of the process just
prior to passage, and the overshoot of a high level. A significant
connection is made with extreme value theory via regular variation or
maximum domain of attraction conditions imposed on the positive tail of the
canonical measure, which are shown to be necessary for the kind of
convergence behaviour we are interested in.
\end{abstract}

% KEYWORDS
% visi is mazosios raides ir pagal abecele
%
\end{frontmatter}

%s1 #&#
\section{Introduction}\label{s0}

The exit time of a L\'evy process $X$ above a horizontal boundary
has been studied extensively in a variety of situations with a view to
relating its distributional behaviour to the tail behaviour of the canonical
measure of $X$. It is helpful to categorise the latter into three general
regimes:
\begin{itemize}
\item {Light tailed} (Cram\'er case).
\item {Medium tailed} (convolution equivalent case).
\item {Heavy tailed} (subexponential tails).
\end{itemize}
 This classification is not prescriptive -- categories may overlap
-- but it provides a convenient general framework in which to summarise
results. Representative papers covering the first two categories are
Bertoin and Doney \cite{BD} for the Cram\'er and Kl\"{u}ppelberg,
Kyprianou and Maller~\cite{kkm} for the convolution equivalent case.
%% (more references are given later).
The intention of the present paper is to consider in some detail the passage
time problem with special emphasis on the third category -- the heavy tailed
cases.

We assume subexponential growth together with regular variation or maximum
domain of attraction conditions for the positive part of the canonical
measure of $X$, or of its increasing ladder height process; on the negative
side, we assume regular variation of the renewal measure of the descending
ladder process, allowing both finite and infinite mean cases. To these is
added the assumption of a drift to $-\infty$ a.s. of the process, possibly
at a linear rate, as is the case when the process has finite mean, but
possibly at a much faster rate.
%%, when the mean is infinite but drift to $-\infty$ still obtains.
We obtain very explicit and detailed descriptions of the asymptotic
behaviours of the process, in these situations. In particular, we obtain
local, and functional, versions of limit distributions for the passage time
itself, as well as for the position of the process just prior to passage,
and the overshoot of a high level.

Our results are original in a number of respects. We give a very general
treatment for L\'evy processes,
%with results phrased in terms of the tail of the canonical measure of
%the process itself or its ladder processes.
%A point of comparison is with the paper of
% Asmussen and Kl\"uppelberg (1996), who deal with ruin event
%calculations, mainly for random walks and the compound Poisson
%process, as used in insurance risk modelling,
% They also consider the case of subexponential tails, but with moment
%and other restrictions which we relax considerably. We treat {\it
%general L\'evy processes}, and
imposing no overt moment conditions, though it will transpire that our
conditions imply the positive tail of the canonical measure is
integrable (a finite mean for the positive jump process).
% We provide \textit{local} as well as \textit{functional} versions of the
%convergence results, so the results are new even in the finite mean
%case. (In the infinite mean case we know only of the paper by Kl
%\"{u}ppelberg and Kyprianou (2006), which deals with a special case.)
Extreme value theory enters via the regular variation or maximum domain of
attraction conditions we impose on the positive tail of the canonical
measure. These are shown to be \textit{necessary} as well as
sufficient for
convergence of the type we investigate. Subsidiary results in Proposition~\ref{nasc} (concerning the convergence of the overshoot for a general
subordinator) and Proposition~\ref{eq} (concerning connections between the
regular variation or maximum domain of attraction behaviour of the upward
ladder height measure as compared with the L\'evy measure of the underlying
process), are also new, and extend the domain of applicability of the paper.

%An important area of application of results like these %%, which we
%can use for guidance,
%is in insurance risk, where positive
%jumps of the process under consideration represent claims on the
%insurance company's assets, while downward trending in the process
%represents premium income.
%%Triggered by discussions on operational risk, following Basel II and
%Solvency II approaches,
%In recent years there has been a recognition
%that operational risk claims in practice may be well modelled by a
%very heavy tailed distribution, perhaps even having an infinite mean
%(see Embrechts and Samorodnitsky (2003), B\"ocker and Kl\"uppelberg
%(2010), and references in both papers).
%On the other hand, so that the company does not face ruin with
%probability 1,
%it is necessary to assume overall drift to $-\infty$ of the process,
%and since it's desirable to place minimum restrictions on income
%growth, we want to allow for the possibility of a heavy tailed
%distribution in the negative direction as well.

% in the L\'evy process case, whereas those of AK are in terms of the
%distribution of the %increments of the random walk, or the
%distribution of the claim size distribution in the %compound Poisson
%case.
% When specialised to their situations, we show in a remark at the end
%of Section~\ref{s2} %that our results contain those of AK.
In the next section, we introduce the setup. The main results are
stated in
Section~\ref{s1a}, and proofs are in Sections~\ref{s3}--\ref{s5}.
The final
Section~\ref{s6} discusses similar results for random walks and compound
Poisson processes.

%
%Our results can be compared with those of Asmussen and Kl\"uppelberg
%(1996), who also deal with the case of subexponential tails, mainly
%for random walks and the compound Poisson process used in insurance
%risk modelling, but with moment and other restrictions which we relax
%considerably. Asmussen and Kl\"uppelberg describe the distributional
%behavior of the process via a portmanteau approach which shows that
%the path of the process up to the time of passage above a high level
%(representing ruin, in the insurance terminology) is, under their
%assumptions, asymptotic to a certain limiting process as measured by
%the total variation distance. From this they deduce limiting
%distributions for the main quantities of interest, including the
%passage time itself, the position of the process just before crossing
%the boundary, and the overshoot of the boundary. Their results are
%proved in detail for the compound Poisson process, with indicated
%extensions to L\'evy processes.
%

%Further comparison of the Asmussen and Kl\"uppelberg results with ours
%is in Section~\ref{s6}.
%%prove the total variation result for random walk,
%%which is extended in Sections \ref{s3} and \ref{s4} to give detailed
%%convergence results for the overshoot, undershoot and passage time
%itself. These are %%generalised in Section~\ref{s5} for L\'evy
%processes.

%%%%%%%%%%%%%%%%%%%%%%%%%%%%%%%%%%%%%%%%%%%%%%%%%%%%%%%%%%%%%%%%%%%%%%%%%%%%%%%%%%%%%%%%%%%%%%%%%%%%%%%%%%%%%%%
%Section~2
%%%%%%%%%%%%%%%%%%%%%%%%%%%%%%%%%%%%%%%%%%%%%%%%%%%%%%%%%%%%%%%%%%%%%%%%%%%%%%%%%%%%%%%%%%%%%%%%%%%%%%%%%%%%%%%

%s2 #&#
\section{Preliminary setting up}\label{s1}

Let $(X_{t})_{t\geq0}$, $X_{0}=0$, be a real-valued L\'{e}vy process
on a
probability space $\{\Omega,\mathcal{F,P\}}$ with triplet $(\gamma
_{X},\sigma_{X}^{2},\Pi_{X})$, where $\gamma_{X}\in\mathbb{R}$,
$\sigma_{X}^{2}\geq0$ and $\Pi_{X}$ is a L\'{e}vy measure on
$\mathbb{R}$.
Throughout, $X$ is assumed to satisfy
%
%e2.1 #&#
\begin{equation}
\lim_{t\rightarrow\infty}X_{t}=-\infty\qquad \mbox{a.s.} \label{L-inf}
\end{equation}
We refer to Bertoin \cite{bert} and Doney \cite{doneystf} for this
notation and the
ensuing notions of fluctuation theory. Denote by $(H_{t})_{t\geq0}$ the
ascending ladder height subordinator generated by $X$. In view of (\ref
{L-inf}), it is defective, obtained from a non-defective subordinator
$\mathcal{H}$ by independent exponential killing with a rate $q>0$
given by $\mathrm{e}^{-q}=P(H_{1}<\infty)$. By this, we mean there is a non-defective
subordinator $\mathcal{H}$ and an independent exponential variable $e_{q}$
with expectation $1/q$ such that $(H_{t})_{0\leq t<L_{\infty}}$ has the
distribution of $(\mathcal{H}_{t})_{0\leq t<e_{q}}$, where $L_{t}$, $t>0$,
is a local time of $X$
%%, viz, the amount of time $X_t$ is equal to $\sup_{0<s\le t}X_t$ up
%till time $t$;
(cf. Bertoin \cite{bert} Lemma~VI.2, page~157). It follows that
%
%e2.2 #&#
\begin{equation}
P ( H_{t}\leq x ) =P ( H_{t}\leq x,t<L_{\infty} )
=\mathrm{e}^{-qt}P ( \mathcal{H}_{t}\leq x ),\qquad t, x>0.
\label{HcH}
\end{equation}

The descending ladder height subordinator, denoted by $(H_{t}^{\ast
})_{t\geq0}$, is the ascending ladder height subordinator
corresponding to
the dual process $(X_t^*)_{t\ge0}: = (-X_t)_{t\ge0}$. Under \eqref{L-inf},
the process $(H_{t}^{\ast})_{t\geq0}$ is proper, and the
corresponding $q^{\ast}=0$.

Let $\Pi_{\mathcal{H}} (\cdot)$ be the L\'{e}vy measure of $\mathcal{H}$,
with tail $\overline{\Pi}_{\mathcal{H}}(x)=\Pi_{\mathcal{H}} \{
(x,\infty)\}$, $x>0$, assumed positive for all $x>0$. Similarly, $\Pi
_{H^{\ast}}(\cdot)$
is the L\'{e}vy measure of ${H^{\ast}}$, with tail $\overline{\Pi}_{H^*}$,
and we write $\mathrm{d}_{\mathcal{H}}$ and $\mathrm{d}_{H^*}$ for
the drift
coefficients of $\mathcal{H}$ and ${H^{\ast}}$. We have $\mathrm
{d}_{\mathcal{H}}=\mathrm{d}_H$ and $\Pi_\cH=\Pi_H$. Let
$\overline{\Pi
}^+_X$ and
$\overline{\Pi}^-_X$ be the positive and negative L\'{e}vy tails of $X$,
equal to $\Pi_X\{(x,\infty)\}$ and $\Pi_X\{(-\infty,-x]\}$, $x>0$.
Write $\Pi_X^{(+)}$ and $\Pi_X^{(-)}$ for $\Pi_X$ restricted to
$(0,\infty)$ and $(-\infty, 0)$, respectively. Assume throughout that
$\overline{\Pi}^+_X(x)>0$
for all $x>0$.

Our results will be phrased in terms of $\Pi_X$, $\Pi_{\mathcal{H}}$,
and $\Pi_{H^{\ast}}$, or, more specifically, in terms of the behaviour
of their
tails for large values. After normalisation, we can regard these as being
the tails of probability distributions.
%%These tails will always be positive for large values in our setups,
%and we will adopt the %%convention when specifying conditions on them,
%we will mean the corresponding
Then %%%, when the denominators in \eqref{distns} below are positive,
a condition applied to the tail of a probability measure can equally be
applied to the tails of the probability measures defined, for example, by
%
%e2.3 #&#
\begin{equation}
\label{distns} \frac{\Pi_X(\mathrm{d} x)\mathbf{1}_{\{x>1\}}}{\overline{\Pi}_X^+(1)} \quad\mbox{and}\quad \frac{\Pi_{\mathcal{H}} (\mathrm{d}
x)\mathbf{1}_{\{x>1\}}}{\overline{\Pi}_{\mathcal{H}}(1)}.
\end{equation}
%
%%,\ {\rm or}\ \frac{\Pi_{H^*} (\rmd x)\mathbf{1}_{\{x>1\}}}{
%\pibar_{H^*}(1)} \ (x\in\R).
We will need certain functionals of these tails, in particular,
%
%e2.4 #&#
\begin{equation}
\label{Adef} A_X^+(x):=\int_{1}^{x}
\overline{\Pi}_X^+(y)\,\mathrm{d} y\quad \mbox {and}\quad A_X^{\ast}(x):=
\int_{1}^{x}\overline{\Pi}_X^-(y)
\,\mathrm{d} y,\qquad x>1
\end{equation}
and
%
%e2.5 #&#
\begin{equation}
\label{Adef2} A_{\mathcal{H}}(x):=\int_0^x
\overline{\Pi}_{\mathcal{H}}(y)\,\mathrm{d} y \quad\mathrm{and}\quad A_{H^*}(x):=
\int_0^x\overline{\Pi }_{H^*}(y)
\,\mathrm{d} y,\qquad x>0,
\end{equation}
which are kinds of truncated or Winsorised means.

Particular classes of tail functions we are interested in are the regularly
varying ones and the class of probability distributions in the maximum
domain of attraction of the Gumbel distribution. Write $\RV(\alpha)$
for the
class of real valued functions regularly varying at $\infty$ with index
$\alpha\in\mathbb{R}$, so that $\RV(0)$ are the slowly varying
functions. We
refer to Bingham, Goldie and Teugels \cite{BGT} for definitions and properties
of regularly varying functions.
% among them being the Potter bounds (Bingham et al. (1987), p.25):
%when $f(x)\in\RV(\alpha)$,
%for all $\eps>0$ there exist $u_0(\eps)$ and $c_->0$, $c_+>0$ such
%that, whenever $u\ge u_0(\eps)$,
%\be\label{potter0}
%c_-\la^{\alpha-\eps}\le\frac{f(\la u)}{f(u)}\le c_+\la^{\alpha+\eps}
% \mbox{ for every } \la>1.
%\ee

Denote the tail of a distribution function $F$ on $[0, \infty)$ by $%
\overline{F}=1-F$, and assume $\overline{F}(u)>0$ for all $u>0$.
$\overline{F}\in\RV(-\beta)$ for some $\beta\in(0,\infty)$ is
equivalent to $F$
being in the \emph{maximum domain of attraction of a Fr\'echet distribution
with parameter $\beta>0$}, denoted $F\in\operatorname{MDA}(\Phi_\beta
)$. A
positive random variable having distribution tail $\overline{F}$
%%, assumed positive for all $u>0$,
is said to be in the \textit{maximum domain of attraction of the Gumbel
distribution}, which we denote as $\operatorname{MDA}(\Lambda)$, with
\textit{auxiliary function} $a(u)>0$, if
%
%e2.6 #&#
\begin{equation}
\label{ZDGumbel} \frac{\overline{F}(u+a(u)x)}{\overline{F}(u)} \to \mathrm{e}^{-x},\qquad x\ge0.
\end{equation}
(Here and throughout, all limits are as $u\to\infty$ unless
otherwise stated.) Useful properties of such distributions can be found
in Bingham, Goldie and Teugels \cite{BGT}, page~410, Resnick \cite
{resnick}, Chapters~0 and 1, Embrechts, Kl\"{u}ppelberg and Mikosch
\cite{EKM}, Chapter~3 and de Haan and Ferreira
\cite{HaanFerr},
Chapter~1.
%%\C{I drop Balkema and de Haan, since they discovered POT properties
%at the same
%%time as Pickands. We make references to both papers.}
%Balkema and de Haan (1974).
In particular, when \eqref{ZDGumbel} holds, $F$ has finite moments of all
orders, and the auxiliary function $a(u)$ satisfies $a(u)=\mathrm{o}(u)$ and is
\textit{self-neglecting}, that is, $a(u+Ka(u))\sim a(u)$ for any fixed $K$.
Typical distributions in $\operatorname{MDA}(\Phi_\beta)$ are the Pareto
distributions, while $\operatorname{MDA}(\Lambda)$ includes the Weibull and
lognormal.

Further, it is well known from extreme value theory [cf.
Theorems~1.1.2, 1.1.3 and 1.1.6 in de Haan and Ferreira \cite{HaanFerr}]
that \eqref{ZDGumbel} can be extended to give that there is a function
$0<a(u)\rightarrow\infty$ and a positive random variable $C$ such that
%
%e2.7 #&#
\begin{equation}
\frac{\overline{F}(u+a(u)x)}{\overline{F}(u)} \to P(C>x),\qquad  x>0, \label{tv5}
\end{equation}
if and only if (for distributions with unbounded support to the right,
as we
have) $F\in\operatorname{MDA}(\Phi_\beta)$ for some $\beta\in
(0,\infty)$,
or $\overline{F}\in$ $\operatorname{MDA}(\Lambda)$. Furthermore, $a(u)$
can be
chosen as $a(u)=u$ in the first case, and as $a(u)= \int_u^\infty
\overline{F}(y)\,\mathrm{d} y/\overline{F}(u)$ (finite) in the second
case, and $C$ has a $\operatorname{Par}(\beta)$ distribution (i.e., a Pareto
distribution with parameter $\beta>0 $) having density $\beta
(1+x)^{-\beta-1}$, $x>0$, in the first
case, and an exponential distribution with unit parameter ($\operatorname{Exp}(1)$) in the
second case.

%See Proposition~\ref{S} below for connecting $\overline{F}(x)\in\MDA(
%\Lambda)$ with
%$\overline{J}\in\MDA(\Lambda)$. Recall that, when it is finite, we
%denote
%$\ov F_I(u)=\int_u^\infty\overline{F}(y)\rmd y$, $u\ge0$.
%Note that
%distributions whose tails are regularly varying are automatically in $%
%\mathcal{S},$ but the corresponding statement for $\Pi_{\mathcal{H}}
%\in\MDA(\Lambda)$ is
%false.

We introduce also the class of \textit{long-tailed distributions},
$\mathcal{L}$, and the subexponential class, $\mathcal{S}$. The
distribution $F$ (or its
tail $\overline{F}=1-F$) is said to be in class $\mathcal{L}$ if
%
%e2.8 #&#
\begin{equation}
\label{classL} \frac{\overline{F}(u+x)}{\overline{F}(u)}\rightarrow1\qquad \mbox{for } x\in (-\infty,\infty),
\end{equation}
while $F$ (or its tail $\overline{F}$) is said to be in the class
$\mathcal{S}$ of \textit{subexponential} distributions if $F\in
\mathcal
{L}$ and
%
%e2.9 #&#
\begin{equation}
\frac{\overline{F^{2\ast}}(u)}{\overline{F}(u)}\rightarrow2, \label{S2}
\end{equation}
where $F^{2\ast}=F\ast F$. For background, see Foss, Korshunov and
Zachary \cite{FKZ}.
%%Kl\"{u}ppelberg et al. (2004), Sect.~1.3.2.
We have $\RV(\alpha)\subset\mathcal{S\subset L}$ but $\operatorname{MDA}(\Lambda)$ is not contained in $\mathcal{S}$ [Goldie and Resnick
\cite{GR}].

Consistent with the convention noted in \eqref{distns}, abbreviate
$\Pi
_X^{(+)} (\mathrm{d} x)\mathbf{1}_{\{x>1\}}/\overline{\Pi
}_X^+(1)\in\operatorname{MDA}(\Lambda)$ to $\Pi_X^{(+)} \in \operatorname{MDA}(\Lambda)$
and $\Pi_\cH(\mathrm{d} x)\mathbf{1}_{\{x>1\}}/\overline{\Pi
}_\cH
(1)\in
\mathcal{S}$ to $\Pi_\cH\in\mathcal{S}$, etc. With this notation, our
second basic assumption is
%
%e2.10 #&#
\begin{equation}
\label{LHsub} \Pi_\cH\in\mathcal{S}.
\end{equation}
Equation~\eqref{LHsub} is equivalent to $P(\mathcal{H}_{1}\in\cdot)\in
\mathcal
{S}$, and then $P(\mathcal{H}_{1}>u)\sim\overline{\Pi}_{\mathcal
{H}}(u)$ as $u\to\infty$
[Embrechts, Goldie and Veraverbeke \cite{EGV}, Pakes \cite
{Pakes,Pakes2}]. Together with \eqref{L-inf}, \eqref{LHsub} implies that
%
%e2.11 #&#
\begin{equation}
\label{max} P\Bigl(\sup_{t\ge0}X_t> u\Bigr) \sim
q^{-1}\overline{\Pi}_{\mathcal{H}}(u)\qquad \mbox{as } u\to\infty
\end{equation}
[from Lemma~3.5 of Kl\"{u}ppelberg, Kyprianou and Maller \cite{kkm}].
%%(Bertoin and Doney (1994)).

For $u>0$ let
%
%e2.12 #&#
\begin{equation}
\tau_{u}:=\inf\{t>0\dvt X_{t}>u\}, \qquad Z^{(u)}=-X_{\tau_u-},\qquad
\mathrm{O}^{(u)}=X_{\tau_u}-u \label{Lptu}
\end{equation}
denote the passage time above level $u>0$, the negative of the position
reached just prior to passage, and the overshoot above the level. (The
reason for taking $-X$ in the definition of $Z$ will become apparent later.)
Note that $P(\tau_{u}<\infty)=P(H_{\infty}>u)<1$ for all $u>0$ by
\eqref{L-inf}, while $P(\tau_{u}<\infty)>0$ for all $u>0$ because of our
assumption that $\overline{\Pi}^+_X(x)>0$ for all $x>0$ and
$\lim_{u\to\infty}P(\tau_u<\infty)=0$ by \eqref{max}. We use
$P^{(u)}(\cdot)=P(\cdot\mid\tau_{u}<\infty)$, $u>0$, defined in an
elementary
way, for the probability measure conditional on passage above $u$. We also
use the notation $\overline{X}_t=\sup_{0<s\le t}X_s$, $t\ge0$.

%Define $G^{\ast}$ as the
%renewal function in the downgoing ladder height process of $X$:
%\be\label{Gdef}
%G^*(x):= \int_0^\infty P(H_t^*\le x)\rmd t, x\ge0.
%\ee
Recall the definition of $A_{H^*}(\cdot)$ in \eqref{Adef2}. Our third
main assumption is of the form:
%
%e2.13 #&#
\begin{equation}
\label{Gcon} A_{H^*}(\cdot)\in\RV(\gamma),
\end{equation}
where the precise value of the index $\gamma\in[0,1)$ will be specified
later. By, for example, Bingham, Goldie and Teugels \cite{BGT},
page~364, \eqref{Gcon} is equivalent to
$G^{\ast}(\cdot)\in\RV(1-\gamma)$, where $G^*$ is the renewal
measure for
the strict decreasing ladder height process, and then we have, as
$x\rightarrow\infty$,
%
%e2.14 #&#
\begin{equation}
\label{m} A_{H^*}(x) \sim\frac{k_{\gamma}x}{G^{\ast}(x)}\in\RV(\gamma )\qquad\mbox{where }k_{\gamma}=\frac{1}{\Gamma(1+\gamma)\Gamma(2-\gamma)}.
\end{equation}
Equation~\eqref{Gcon} is also equivalent to
%
%e2.15 #&#
\begin{equation}
\label{xhA} \lim_{x\to\infty}\frac{x\overline{\Pi
}_{H^*}(x)}{A_{H^*}(x)}=\gamma,\qquad 0\le
\gamma<1
\end{equation}
(Bingham, Goldie and Teugels \cite{BGT}, Theorem~1.5.11, page~18,
Theorem~1.6.1, page~30).
% \eqref{Gcon} is equivalent to
%\begin{equation}
%A_{H^*}(x)\sim\frac{k_{\gamma}x}{G^{\ast}(x)}\in\RV(\gamma),\mbox{
%where }k_{\gamma}=\frac{1}{\Gamma(1+\gamma)\Gamma(2-\gamma)}.
%\label{tv22}
%\end{equation}%
%(??ref??)

%Our main objective is to prove results
%about the asymptotic behaviour of $X$ under $P^{(u)}$ as $u\rightarrow
%\infty$ in some generality.

%s3 #&#
\section{Main results}\label{s1a}

We now state our two main results. Both assume (\ref{L-inf}) and
(\ref{LHsub}), and the first assumes in addition that $A_{H^*}\in\RV(0)$,
that is, that $A_{H^*}$ is slowly varying as $x\to\infty$. This
implies that
$X^{\ast}_t$ is positively relatively stable as $t\to\infty$, so there
is a continuous, strictly increasing function $c(\cdot)\in\RV(1)$
such that $X_{t}^{\ast
}/c(t)\overset{\mathrm{P}}{\longrightarrow}1$ as $t\rightarrow
\infty$.
This in turn implies that the process $ ( X^{\ast
}_{st}/c(t)
)_{0\le s\le1}$ converges weakly in $\mathbb{D}_{0}[0,1]$ (i.e., in
the sense
of weak convergence of c\`adl\`ag functions on $[0,1]$ with the
Skorokhod topology) as $t\to\infty$ to the process $\mathbf{D}^{(0)}$,
where $\mathbf{D}^{(0)}(s)\equiv s$. This situation includes the
possibility of a finite,
positive mean for $X_1^*$.
%generalising a corresponding case in \cite{AK}.
Write $b(\cdot)$ for the inverse function of $c(\cdot)$.
We sometimes write $X^*(t)$ for $X^*_t$.

%th3.1 #&#
\begin{thmm}\label{old} Assume
$\lim_{t\to\infty}X_{t}=-\infty$ a.s., $\Pi_\cH\in\mathcal{S}$, and
$A_{H^*}\in\RV(\gamma)$ with $\gamma=0$.

%% is slowly
%varying
%as $x\to\infty$.
(1)   Then the following are equivalent;
\begin{longlist}[(a)]
\item[(a)]
there exists $a(u)>0$ with $\lim_{u\to\infty}a(u)=\infty$ such that
$P^{(u)}(\mathrm{O}^{(u)}\in a(u)\,\mathrm{d}x)$, $x>0$, has a non-degenerate
limit as $u\rightarrow\infty$;
\item[(b)]
either $\overline{\Pi}_{\mathcal{H}}\in\RV(1-\gamma-\beta)$ for some
$\beta>1-\gamma$
and then \textup{(a)} holds with $a(u)=u$ [case~\textup{(i)}] or else $\Pi_{\mathcal{H}}
\in\operatorname{MDA}(\Lambda)$, and then \textup{(a)} holds with $a(u)=\int_u^\infty
\overline{\Pi}_\cH(y)\,\mathrm{d} y/\overline{\Pi}_\cH(u)$ [case \textup{(ii)}];
\item[(c)]
either $\overline{\Pi}^+_X\in\RV(-\beta)$ for some $\beta>1$ (case
\textup{(i)}) or else $\Pi_X^{(+)} \in\operatorname{MDA}(\Lambda)$ [case \textup{(ii)}],
and $a(\cdot)$ may then be chosen as $a(u)=u$ in the first case or as
$a(u)=\int_u^\infty\overline{\Pi}_X^+(y)\,\mathrm{d} y/\overline
{\Pi}_X^+(u)$
in the second case.
\end{longlist}
(2)   When \textup{(a)}--\textup{(c)} hold, the $P^{(u)}$-distribution of $\tau_u$,
restricted to
the event $X_{\tau_u-}<u$, has a density $g^{(u)}(\cdot)$ which satisfies
%%% with $r(u)=b(a(u)),$%
%e3.1 #&#
\begin{equation}\label{local}
\lim_{u\rightarrow\infty}b\bigl(a(u)\bigr)g^{(u)}\bigl(tb\bigl(a(u)
\bigr)\bigr)= \cases{ %
\displaystyle\frac{\beta-1}{(1+t)^{\beta}}, & \quad$\mbox{in
case \textup{(i)},}$
\vspace*{2pt}\cr
\mathrm{e}^{-t}, &\quad $\mbox{in case \textup{(ii)},}$}
\end{equation}
uniformly on compact subintervals of $(0,\infty)$. Moreover,
conditioned on $\tau_u=tb(a(u))$, the $P^{(u)} $-finite-dimensional
distributions of the process
\[
\biggl\{ \frac{X^{\ast}(s\tau_u)}{c(\tau_u)}, 0\le s\le1 \biggr\}
\]
converge to those of $\mathbf{D}^{(0)}$ as $u\to\infty$.

(3)   Further: when \textup{(a)--(c)} hold, under $P^{(u)}$ the process
%
%e3.2 #&#
\begin{equation}
\mathbf{Y}^{(u)}:= \biggl( \frac{Z^{(u)}}{a(u)}, \frac
{\mathrm{O}^{(u)}}{a(u)},
\frac{\tau_u}{b(a(u))}, \biggl( \frac{X^{\ast}(s\tau
_u)}{a(u)} \biggr) _{0\le s\le1} \biggr)
\label{Y}
\end{equation}
converges weakly as $u\to\infty$ in $\mathbb{R}^3\times\mathbb
{D}_{0}[0,1]$ to $ (
V,U,V,(V\mathbf{D}^{(0)}(s))_{0\le s\le1} )$, where
in case \textup{(i)}
%
%e3.3 #&#
\begin{equation}
\label{Y0} P (V\in\mathrm{d}z,U\in\mathrm{d}x )=\frac{\beta(\beta
-1)\,\mathrm{d}z\,\mathrm{d}x}{(1+z+x)^{\beta+1}}, \qquad x,z>0,
\end{equation}
and in case \textup{(ii)}
%
%e3.4 #&#
\begin{equation}
\label{Y1} P (V\in\mathrm{d}z,U\in\mathrm{d}x )=\mathrm{e}^{-z-x}\,\mathrm {d}z
\,\mathrm{d}x,\qquad  x,z>0.
\end{equation}
\end{thmm}

%re3.1 #&#
\begin{rem}\label{3.1} (i)
The redundant parameter $\gamma=0$ is introduced in Theorem~\ref{old}
for conformity with Theorem~\ref{new}, below.

\textup{(ii)}
The event $\{X_{\tau_u-}<u\}$ in Theorem~\ref{old} has
$P^{(u)}$-probability approaching 1 as $u\to\infty$; see Remark~\ref
{ont} in Section~\ref{s4}.

\textup{(iii)} In general, we cannot replace condition \eqref{LHsub}
with simpler equivalent conditions on $\Pi_X$ directly, but easily checked
sufficient conditions are available; see Remark~\ref{6.1} in Section~\ref{s5}.

(iv) The assumption $A_{H^*}\in\RV(0)$ in Theorem~\ref{old} is true
in particular when $0<A_{H^*}(\infty)<\infty$, or, equivalently, when
$0<EX_1^*<\infty$, so the case of a finite mean for $EX_1^*$ is
included in the theorem. Note that part 1\textup{(c)} implies $EX_1^+=E(X_1\vee
0)<\infty$ in any case.
A related result for random walks and compound Poisson processes
with finite mean is in Asmussen and Kl\"uppelberg~\cite{AK}.
\end{rem}

In our next result, we replace the assumption $A_{H^{\ast}}\in\RV(0)$
by the
condition that $A_{H^{\ast}}\in\RV(\gamma)$ for some $\gamma\in(0,1)$.
This can only happen when $E|X_{1}|=\infty$, and we will show that it
is in
fact equivalent, under our basic assumptions,\vadjust{\goodbreak} to $\overline{\Pi
}_{X}^{-}\in
\RV(\gamma-1)$ [see Proposition~\ref{Q}, where $A_{H^{\ast}}$ is
shown to
be asymptotically equivalent to $q^{-1}A_{X}^{\ast}$, and note \eqref
{xpA}]. It then follows that $X^{\ast}$ is in the domain of attraction
of $\mathbf{D}$, a standard stable
subordinator of parameter $\overline{\gamma}:=1-\gamma\in(0,1)$. Let
$c(\cdot)$ be such that $ ( X_{st}^{\ast}/c(t) ) _{0\le
s\le
1}\overset{D}{\rightarrow}\mathbf{D}$ as $t\to\infty$, and let
$b(\cdot)$ denote the inverse function of\vspace*{1pt}
$c(\cdot)$, so that $b(\cdot)\in\RV(\overline{\gamma})$, and let
$\widehat{\mathbf{D}}_{t},_{z}$ denote an associated
\textquotedblleft stable
subordinator bridge'', which is a rescaled version of $\mathbf{D}$
conditioned to be at $z>0$ at time $t$; namely,
\[
P ( \widehat{\mathbf{D}}_{t},_{z}\in\mathcal{B} ) =P
\bigl( \bigl(D(ts)\bigr)_{0\le s\le1}\in\mathcal{B}\mid D_{t}=z
\bigr),
\]
for any Borel set $\mathcal{B}$. Thus, with
%
%e3.5 #&#
\begin{equation}
h_{t}(x)\,\mathrm{d}x=P(D_{t}\in\mathrm{d}x) \label{hdef}
\end{equation}
as the density of $D$, we have for $0=s_{0}<s_{1}<s_{2}<\cdots
<s_{k}<1$, $y_{0}=0$, and $y_{1}<y_{2}<\cdots<y_{k}<z$,
%
%e3.6 #&#
\begin{equation}
P \Biggl( \bigcap_{r=1}^{k}\bigl\{
\widehat{D}_{t,z}(s_{r})\in\mathrm {d}y_{r}\bigr
\} \Biggr) =\frac{h_{t(1-s_{k})}(z-y_{k})}{h_{t}(z)}\prod_{r=1}^{k}h_{t(s_{r}-s_{r-1})}(y_{r}-y_{r-1})
\,\mathrm{d}y_{r}. \label{12}
\end{equation}
We will use $\widehat{\mathbf{D}}_{W,V}$ in the obvious sense, where $(W,V)$
are positive random variables independent of the family $\widehat
{\mathbf{D}}_{t},_{z}$.

%th3.2 #&#
\begin{thmm}
\label{new} Assume $\lim_{t\to\infty}X_{t}=-\infty$ a.s., $\Pi
_\cH\in
\mathcal{S}$, and $A_{H^*}\in\RV(\gamma)$ with $\gamma\in(0,1)$.

\begin{longlist}[(1)]
\item[(1)]   Then conditions \textup{(a)--(c)} of Theorem~\ref{old} remain equivalent
as stated for the current value of $\gamma\in(0,1)$.

%(a) $P^{(u)}(\mathrm{O}^{(u)}\in a(u)\mathrm{d} x)$ has a non-degenerate limit
%for
%some $a(u)>0$, $a(u)\rightarrow\infty$;
%
%(b) either $\overline{\Pi}_{\mathcal{H}}\in\RV(1-\gamma-\beta)$ for
%some $%
%\beta>1-\gamma$ and then (a) holds with $a(u)=u$ (case (i)), or else
%$\Pi_{%
%\mathcal{H}} \in MDA(\Lambda)$ and then (a) holds with $a(u)=\int_u^
%\infty
%\overline{\Pi}_\cH(y) \mathrm{d} y/\overline{\Pi}_{\mathcal{H}}(u)$
%(Case\textup{(ii)});
%
%(c) either $\overline{\Pi}^+_X\in\RV(-\beta)$ for some $\beta>1-
%\gamma$
%(Case (i)) or else $\Pi_X^{(+)}\in\operatorname{MDA}(\Lambda)$ (Case
%\textup{(ii)}),
%and $a(u)$ may then be chosen as $a(u)=u$ in the first case or as $%
%a(u)=\int_u^\infty\overline{\Pi}_X^+(y)\mathrm{d} y/\overline{
%\Pi}_X^+(u)$
%in the second case.
\item[(2)]   Assume conditions \textup{(a)--(c)} as stated in Theorem~\ref{old} hold
for the current value of $\gamma\in(0,1)$, and further assume that
$X_{t}$ has a non-lattice distribution for each fixed $t>0$. Then,
uniformly for $z\in{}[\Delta
_{0},\Delta_{1}]$, for any fixed $0<\Delta_{0}<\Delta_{1}<\infty$,
and $t\in{}[ T_{0},T_{1}]$ for any fixed $0<T_{0}<T_{1}<\infty$,
%
%e3.7 #&#
\begin{equation}
\lim_{u\rightarrow\infty}a(u)b\bigl(a(u)\bigr)P^{(u)}\bigl(Z^{(u)}
\in \bigl(a(u)z,a(u)z+\Delta \bigr],\tau_{u}\in b\bigl(a(u)\bigr)\,\mathrm{d}t
\bigr)=h_{t}(z)f(z)\Delta\,\mathrm{d}t, \label{Z}
\end{equation}
where, in case \textup{(i)},
%
%e3.8 #&#
\begin{equation}
f(z)=\frac{\Gamma(\beta)}{\Gamma(\beta+\gamma-1)(1+z)^{\beta}}, \label{fdef}
\end{equation}
and in case \textup{(ii)}
\[
f(z)=\mathrm{e}^{-z},\qquad z>0.
\]
Moreover, for $k=2,3,\ldots,$ take $z_i>0$ and
$I_{i}=(a(u)z_{i},a(u)z_{i}+\Delta_{i}]$, $i=1,2,\ldots,k-1$, and
write, for $0<s_{1}<\cdots<s_{k-1}<s_k=1$,
\[
A_{k}= \bigl\{ X^{\ast}\bigl(s_{i}tb\bigl(a(u)
\bigr)\bigr)\in I_{i}, i=1,2,\ldots ,k-1 \bigr\}.
\]
Then, uniformly for %%$\Delta_{i}\in(0,\Delta_{0}]$, and
$z_{i}\in{}[\Delta_{0},\Delta_{1}]$, $i=1,2,\ldots,k$, for any
fixed $0<\Delta_{0}<\Delta_{1}<\infty$, and $t\in{}[ T_{0},T_{1}]$
for any fixed $0<T_{0}<T_{1}<\infty$, we have
%
%e3.9 #&#
\begin{eqnarray}\label{theta}
&&\lim_{u\rightarrow\infty}\bigl(a(u)\bigr)^{k}b\bigl(a(u)
\bigr)P^{(u)}\bigl(A_{k},Z^{(u)}\in \bigl(z_ka(u),
z_ka(u)+\Delta_{k}\bigr], \tau_{u}\in b\bigl(a(u)
\bigr)\,\mathrm{d}t\bigr)
\nonumber
\\[-8pt]
\\[-8pt]
\nonumber
&&\quad =\theta(z_{1},z_{2},\ldots,z_{k},t)\prod
_{i=1}^{k}\Delta _{i}
\,\mathrm{d}t,\qquad k=1,2,\ldots.
\end{eqnarray}
Here, with $s_0=z_{0}=0$,
\[
\theta(z_{1},z_{2},\ldots,z_{k},t)=\prod
_{i=1}^{k}h_{t(s_{i}-s_{i-1})}(z_{i}-z_{i-1})f(z_{k}).
\]
\item[(3)]   Further: assume conditions \textup{(a)--(c)} as stated in Theorem~\ref
{old} hold for the current value of $\gamma\in(0,1)$, and that $X_{t}$
has a non-lattice distribution
for each $t>0$. Then, under $P^{(u)}$, the process $\mathbf{Y}^{(u)}$
defined in (\ref{Y}) converges weakly in $\mathbb{R}^3\times\mathbb
{D}_{0}[0,1]$ as $u\to\infty$ to the process $ (V,U,W,(\widehat
{\mathbf{D}}_{W,V}(s))_{0\le s\le1} ) $, where
in case \textup{(i)}
%
%e3.10 #&#
\begin{eqnarray}
\label{J11}&& P (V\in\mathrm{d}z,U\in\mathrm{d}x,W\in\mathrm{d}t )
\nonumber
\\[-8pt]
\\[-8pt]
\nonumber
&&\quad=
\frac
{\Gamma(\beta+1)}{\Gamma(\beta+\gamma-1)(1+z+x)^{\beta
+1}}h_{t}(z)\,\mathrm{d}z \,\mathrm{d}x \,\mathrm{d}t,\qquad t,x,z>0,
\end{eqnarray}
and in case \textup{(ii)}
%
%e3.11 #&#
\begin{equation}
P (V\in\mathrm{d}z,U\in\mathrm{d}x,W\in\mathrm{d}t )=\mathrm{e}^{-z-x}h_{t}(z)
\,\mathrm{d}z \,\mathrm{d}x \,\mathrm{d}t,\qquad  t,x,z>0. \label{J12}
\end{equation}
\end{longlist}
\end{thmm}

%re3.2 #&#
\begin{rem}
\textup{(i)} The further assumption in part 2 of Theorem~\ref{new}, that for
each $t>0$, $X_{t}$ has a non-lattice distribution, is equivalent to
assuming that $X$ is not a compound Poisson process whose step distribution
takes values on a lattice. We can cover the lattice case also with only
minor adjustments. Thus, if the lattice has span 1, we need only
restrict $\Delta$ to take integer values and replace
$(a(u)z,a(u)z+\Delta]$ in (\ref{Z}) by $(\lfloor a(u)z\rfloor
,\lfloor
a(u)z\rfloor+\Delta]$, and similarly in (\ref{theta}), for a valid
conclusion. The only difference in the proof is
which version of a local limit theorem is used.

\textup{(ii)} The right-hand sides of \eqref{Y0} and \eqref{Y1} and \eqref{J11}
and \eqref{J12} are probability densities on $x,z>0$ and $t,x,z>0$, so,
under the conditions of Theorems \ref{old} and \ref{new},
the limiting distributions of $Z^{(u)}/a(u)$ (and of course those of
$\mathrm{O}^{(u)}/a(u)$ and $\tau_u/b(a(u))$) are concentrated on $[0,\infty)$.
Thus, $\lim_{u\to\infty}P(Z^{(u)}/a(u)\le-z)=0$ for all $z>0$. So it
is convenient
to define $Z^{(u)}=-X_{\tau_u-}$ as we did in \eqref{Lptu}.

%\textrm{\textup{(iii)}\ Since $\lim_{u\to\infty}P(Z^{(u)}/a(u)\le0)=0$ in
%Theorems \ref{old} and \ref{new}, $\lim_{u\to\infty}P(X_{\tau_u}=u,
%\tau_u<\infty)=0$, that is, $X$ creeps over level $u$
%with probability tending to 0 as $u\to\infty$. This follows because, in
%order to creep with $Z^{(u)}>0$, $X$ would have to pass continuously
%over
%the interval $(0,u)$, or, equivalently $\mathcal{H}$ would have to
%reach
%level $u$ without any jumps. This probability is exponentially small,
%or
%zero if $(0,\infty)$ is regular for 0. It follows then from Lemma~5.1
%of Griffin and Maller (2011) that
% $\lim_{u\to\infty}P(X_{\tau_u-}=u,X_{\tau_u}>u, \tau_u<\infty)=0$,
%hence also $\lim_{u\to\infty}P(X_{\tau_u-}=u, \tau_u<\infty)=0$. So $
%\lim_{u\to\infty}P(X_{\tau_u-}<u, \tau_u<\infty)=1$ and the
%restriction to $X_{\tau_u-}<u$ is redundant.
%}

\textup{(iii)} In connection with Theorem~\ref{new}, we mention the paper by
Kl\"{u}ppelberg and Kyprianou~\cite{kk}, which deals with the infinite
mean case
under special assumptions.

(iv) The marginal limiting distributions of the fluctuation
quantities are easily computed from \eqref{Y0} and \eqref{Y1} and
\eqref{J11} and \eqref{J12}. The identities $t^{1/\overline{\gamma}
}h_t(z)=h_1(z/t^{1/\overline{\gamma}})$ and $\int_0^\infty h_t(z)\,\mathrm{d}
t=z^{-\gamma}/\Gamma(\overline{\gamma})$, where $\overline{\gamma
}=1-\gamma$
[see Sato \cite{sato}, pages~87, 261)], are useful. Thus, for example,
under the
conditions of case~\textup{(i)} of Theorem~\ref{new}, the limiting densities of
$(Z^{(u)}, \mathrm{O}^{(u)})$ and $\tau_u$, suitably normalised, are derived
from \eqref{J11} as
%
%e3.12 #&#
\begin{equation}
\label{J111} P (V\in\mathrm{d}z,U\in\mathrm{d}x )= \frac{\Gamma
(\beta
+1)z^{-\gamma}}{\Gamma(1-\gamma)\Gamma(\beta+\gamma
-1)(1+z+x)^{\beta
+1}}
\,\mathrm{d} z \,\mathrm{d}x,\qquad y,x>0
\end{equation}
and
%
%e3.13 #&#
\begin{equation}
\label{J112} P (W\in\mathrm{d}t )= \frac{\Gamma(\beta)}{\Gamma
(\beta+\gamma
-1)}\int_0^\infty
\frac{h_1(z) \,\mathrm{d}z}{(1+t^{1/\overline
{\gamma
}}z)^\beta} \,\mathrm{d} t,\qquad t>0.
\end{equation}
It can be checked that no pair of $(V,U,W)$ are independent, in case \textup{(i)}.
For case \textup{(ii)},
%
%e3.14 #&#
\begin{equation}
\label{J113} P (V\in\mathrm{d}z,U\in\mathrm{d}x )= \frac{z^{-\gamma
}\mathrm{e}^{-z-x}}{\Gamma(1-\gamma)}
\,\mathrm{d}z \,\mathrm{d}x, \qquad x,z>0
\end{equation}
and
%
%e3.15 #&#
\begin{equation}
\label{J114} P (W\in\mathrm{d}t ) =\int_0^\infty
\mathrm{e}^{-zt^{1/\overline
{\gamma
}}}h_1(z)\,\mathrm{d} z \,\mathrm{d} t,\qquad t>0.
\end{equation}
In this case, $V$ is independent of $U$, $U$ is independent of $W$, but $V$
is not independent of $W$.
\end{rem}

%s4 #&#
\section{Preliminaries to the proofs}\label{s3}

%%We start with a result about $H,$ which actually holds for
%%It is of course
%%relevant because the overshoot of $X$ over level $u>0$ coincides with
%the
%%overshoot of its ladder-height subordinator over the same level.
%% satisfies (\ref{LHsub}).
Our first proposition applies to any defective subordinator, so we change
notation slightly just for this result.

%pr4.1 #&#
\begin{proposition}
\label{nasc} Let $Y$ be any defective subordinator, obtained from a
non-defective subordinator $\mathcal{Y}$ with killing rate $q$, whose
L\'{e}vy measure is $\Pi_Y$, with tail $\overline{\Pi}_Y$. Assume
$\Pi
_Y\in
\mathcal{S}$. Write $P_Y^{(u)}$ for $P(\cdot\mid T_u^Y<\infty)$,
where $T_u^Y=\inf\{t\dvt Y_{t}>u\}$, $u>0$, and put $\mathrm{O}_Y^{(u)}=Y_{T_u^Y}-u$
on the
event $\{T_u^Y<\infty\}$.

Then $P_Y^{(u)}(\mathrm{O}_Y^{(u)}\in a(u)\,\mathrm{d} x)$
has a non-degenerate limit $P(\mathrm{O}\in\,\mathrm{d} x)$ for some $a(u)>0$,
$a(u)\rightarrow\infty$, if and only if either $\overline{\Pi}_Y\in
\RV(-\alpha)$ for
some $\alpha>0$, or $\Pi_Y \in\operatorname{MDA}(\Lambda)$.

Moreover, in
the first case we can take $a(u)=u$ and $\mathrm{O}$ to have density $\alpha
(1+x)^{-1-\alpha}$, and in the second case we can take $a(u)=\int_u^\infty
\overline{\Pi}_Y(y)\,\mathrm{d} y/\overline{\Pi}_Y(u)=\mathrm{o}(u)$ and $\mathrm{O}$
to have
density $\mathrm{e}^{-x}$.
\end{proposition}

\begin{pf}
For the distribution of $\mathrm{O}_Y^{(u)}$,
%% decompose according to $\{T_u^Y =t\}$ and
use of the compensation formula for Poisson point processes as in
Bertoin \cite{bert}, Proposition~2, page~76, or Kl\"{u}ppelberg,
Kyprianou and Maller \cite{kkm}, Theorem~2.4, gives
\begin{eqnarray*}
P\bigl(\mathrm{O}_Y^{(u)}>xa(u),T^Y_u<\infty
\bigr)&=&P\bigl(Y_{T_u^Y}>u+xa(u),T_u^Y<\infty
\bigr)
\\
&=&E\sum_{0<t<L_\infty}\mathbf{1}_{\{Y_t>u+xa(u),T_u^Y =t\}}
\\
&=&\int_0^\infty \mathrm{e}^{-qt}\int
_{(0,u]} \overline{\Pi }_Y\bigl(u+xa(u)-y\bigr)P (
\mathcal{Y}_{t}\in\mathrm{d} y )\,\mathrm{d} t.
\end{eqnarray*}
From this, writing $e(q)$ for an independent $\operatorname{Exp}(q)$ random variable, we
have for any $C_0>0$
%
%e4.1 #&#
\begin{eqnarray}
\label{Yp}&& P\bigl(\mathrm{O}_Y^{(u)}>xa(u), T_u^Y<
\infty\bigr)\nonumber\\
&&\qquad = q^{-1}\int_{(0,u]}P(
\mathcal{Y}_{e(q)}\in\mathrm{d} y)\overline{\Pi }_Y
\bigl(u+xa(u)-y\bigr)
\\
&&\qquad= q^{-1} \biggl( \int_{(0,C_0]}+\int
_{(C_0,u]} \biggr) P(\mathcal {Y}_{e(q)}\in\mathrm{d} y)
\overline{\Pi}_Y\bigl(u+xa(u)-y\bigr).\nonumber
\end{eqnarray}
Assume at this stage that $\Pi_\cY\in\mathcal{S }$. Then $\Pi_\cY
\in
\mathcal{L}$, so we have
%
%e4.2 #&#
\begin{equation}
\label{sxz} \overline{\Pi}_\cY\bigl(u-y+xa(u)\bigr)\sim\overline{
\Pi}_\cY\bigl(u+xa(u)\bigr)\qquad \mbox{uniformly for } y
\in(0,C_{0}] \mbox{ and } x \ge0.
\end{equation}
Thus,
%
%e4.3 #&#
\begin{equation}
\label{Hin} \int_{(0,C_0]}P(\mathcal{Y}_{e(q)}\in
\mathrm{d} y)\overline{\Pi }_\cY \bigl(u+xa(u)-y\bigr) \sim P(
\mathcal{Y}_{e(q)}\le C_{0})\overline{\Pi}_\cY
\bigl(u+xa(u)\bigr).
\end{equation}
%
%\nonumber\\
%&\sim&
%P(\cY_{e(\varkappa)}\le C_{0})P(\mathrm{O}>x)\pibar_\cY(u),
%%\eea
%where the last equivalence follows by \eqref{pirat}.
%Clearly
%\[
%\int_{C_{0}}^{u}P(\cY_{e(q)}\in\rmd y)\pibar_\cY(u+xa(u)-y)
%\leq
%\int_{C_{0}}^{u}P(\cY_{e(q)}\in\rmd y)\pibar_\cY(u-y)
%=qP(T_u^Y<\infty) \sim\pibar_Y(u),
%\]%
Since $\Pi_\cY\in\mathcal{S }$, we know from Lemma~3.5 of Kl\"
{u}ppelberg, Kyprianou and Maller \cite{kkm} (with $\alpha=0$) that
$\overline{\Pi}_\cY(u)\sim qP(T_u^Y<\infty)$. Given arbitrary
$\varepsilon\in(0,1)$, we can choose $C_0>0$ such that $P(\mathcal
{Y}_{e(q)}>C_{0})\leq\varepsilon$. Then for $u$ large enough,
again using \eqref{sxz}, %%the fact that $\Pi_\cY\in\cL$,
\begin{eqnarray*}
(1+\varepsilon) \overline{\Pi}_\cY(u) &\ge& qP\bigl(T_u^Y<
\infty\bigr)
\\
&=& \biggl(\int_{(0,C_0]}+\int_{(C_{0},\infty)} \biggr)
P(\mathcal {Y}_{e(q)}\in\mathrm{d} y)\overline{\Pi}_\cY(u-y)
\nonumber
\\
&\ge& (1-\varepsilon)P(\mathcal{Y}_{e(q)}\le C_{0})
\overline{\Pi }_\cY (u)\\
&&{}+\int_{(C_0,u]}P(
\mathcal{Y}_{e(q)}\in\mathrm{d} y)\overline {\Pi }_\cY(u-y),
\end{eqnarray*}
giving
\[
\int_{(C_0,u]}P(\mathcal{Y}_{e(q)}\in\mathrm{d} y)
\overline{\Pi }_\cY(u-y) \le\bigl((1+\varepsilon)-(1-
\varepsilon)^2\bigr) \overline{\Pi}_\cY(u)\le 3\varepsilon
\overline{\Pi}_\cY(u).
\]
From this, and \eqref{Yp} and \eqref{Hin}, and since $\overline{\Pi
}_\cY
(u)\sim qP(T_u^Y<\infty)$, we have
%
%e4.4 #&#
\begin{eqnarray}
\label{Yp2} P^{(u)} \bigl(\mathrm{O}_Y^{(u)}>xa(u)
\bigr)&=& \frac{P (\mathrm{O}_Y^{(u)}>xa(u),
T_u^Y<\infty )}{P(T_u^Y<\infty)}
\nonumber
\\[-8pt]
\\[-8pt]
\nonumber
&=& \bigl(1+\mathrm{o}(1)\bigr) P(\mathcal{Y}_{e(q)}\le C_{0})
\frac{\overline{\Pi
}_Y(u+xa(u))}{\overline{\Pi}_Y(u)} +\mathrm{o}(1).
\end{eqnarray}
As discussed in \eqref{tv5}, the condition $\overline{\Pi}_Y\in\RV
(-\alpha)$
for some $\alpha>0$, or $\Pi_Y \in\operatorname{MDA}(\Lambda)$, is
equivalent to the existence of $a(u)\rightarrow\infty$ such that
%
%e4.5 #&#
\begin{equation}
\label{pirat} \frac{\overline{\Pi}_Y(u+xa(u))}{\overline{\Pi}_Y(u)}\rightarrow P(\mathrm{O}>x),
\end{equation}
and when it holds $a(u)$ and $\mathrm{O}$ have the stated properties. The conclusions
of the proposition then follow from this and \eqref{Yp2}.
\end{pf}

We will make use of the ``\'{e}quations amicales'' of Vigon \cite{Vig}, which
are
%
%e4.6 #&#
\begin{equation}
\overline{\Pi}^+_X(u)=\int_{(0,\infty)}\overline{
\Pi}_{H^*}(y) \Pi _{\mathcal{H}}(u+\mathrm{d}y)+\mathrm{d}_{H^*}n(u),\qquad
u>0 \label{Vig2}
\end{equation}
and
%
%e4.7 #&#
\begin{equation}
\overline{\Pi}^-_X(u)=\int_{(0,\infty)}\overline{
\Pi}_{\mathcal{H}}(y) \Pi_{H^{\ast}}(u+\mathrm{d}y)+\mathrm{d}_{\mathcal{H}}n^{\ast
}(u)+q
\overline{\Pi}_{H^*}(u),\qquad u>0, \label{Vig3}
\end{equation}
%
%%which can, equivalently, be written as
%%\be\label{Vig4}
%%\pX(u)=\int_{(0,\infty)}\left(\pcH(u)-\pcH(y+u)\right)\Pi_{\starH}(
%\rmd y)+{\rm
%%d}_{\starH}\Pi'_{\cal H}(u),\ee(compare \eqref{WHf}),
where $n(\cdot)$, $n^{\ast}(\cdot)$ denote c\`adl\`ag versions of the
densities of $\Pi_{\mathcal{H}}$, $\Pi_{H^{\ast}}$, defined if
$\mathrm
{d}_{\mathcal{H}}>0$, $\mathrm{d}_{H^*}>0$, respectively.
Recall that $q$ is the killing rate in \eqref{HcH}.

We are looking for limit theorems which will always include the convergence
of the normed overshoot,
%%so from now on we will add to our basic assumptions the following:%
and Proposition~\ref{nasc} suggests the relevance of conditions like
%
%e4.8 #&#
\begin{equation}
\overline{\Pi}_{\mathcal{H}} \in\RV(-\alpha)\qquad\mbox{for some }\alpha>0
\mbox{ [case (i)]}\quad \mbox{or}\quad  \Pi_{\mathcal{H}}\in\operatorname{MDA}(\Lambda)\qquad \mbox{[case \textup{(ii)}]}. \label{HS}
\end{equation}
%
%Results in Asmussen and Kl\"{u}ppelberg (1996) suggest that when
%$E|X_{1}|<\infty$,
%so that $EX_{1}\in(-\infty,0)$, and $EH_{1}^{\ast}<\infty$,
%we have (\ref{HS}) equivalent to
%\begin{equation} \label{mf}
%\pibar^+_X\mbox{ }\in\RV(-\beta)\mbox{ for some }\beta>1\mbox{ or }%
%\Pi_X \in\MDA(\Lambda).
%\end{equation}%
%We will prove this, and in fact a more general result, in t
The next proposition shows that these imply similarly stated conditions
on ${\Pi}_X^{(+)}$. At this stage, we are not assuming $\Pi_\cH\in
\mathcal{S}
$.

%pr4.2 #&#
\begin{proposition}\label{eq} Assume $\lim_{t\to\infty
}X_{t}=-\infty$
a.s. and $A_{H^*}\in
\RV(\gamma)$ with $\gamma\in{}[0,1)$. Suppose (\ref{HS}) holds
with $\alpha=\beta+\gamma-1>0$, where $\beta>0$, in case \textup{(i)}.

Then $\overline{\Pi}^+_X$ $\in
\RV(-\beta)$ (case \textup{(i)}), or $\Pi_X^{(+)} \in\operatorname{MDA}(\Lambda
)$ [case \textup{(ii)}], or, equivalently,
%
%e4.9 #&#
\begin{equation}
\label{pX} \frac{\overline{\Pi}_X^+(u+xa(u))}{\overline{\Pi
}_X^+(u)}\rightarrow P(C>x),\qquad x>0,
\end{equation}
where $a(u)=u$ and $P(C>x)=(1+x)^{-\beta}$ (case \textup{(i)}), or $a(u)=\int_u^\infty\overline{\Pi}_\cH(y)\,\mathrm{d} y/\overline{\Pi
}_{\mathcal
{H}}(u)$ and $P(C>x)=\mathrm{e}^{-x}$ [case \textup{(ii)}]. Further, in both cases we
have, for
some constants $c_{\gamma,\beta}\in(0,\infty)$ (whose values are made
explicit in the proof),
%
%e4.10 #&#
\begin{equation}
\overline{\Pi}^+_X(u)\sim\frac{c_{\gamma,\beta}\overline{\Pi
}_{\mathcal{H}}(u)A_{H^*}(a(u))}{a(u)}. \label{est}
\end{equation}
Moreover, in case \textup{(ii)} we can alternatively take $a(u)=\int_u^\infty
\overline{\Pi}^+_X(y)\,\mathrm{d} y/\overline{\Pi}^+_X(u)$, $u>0$.
\end{proposition}

%\begin{rem}
%It is actually true that under (\ref{L-inf}), (\ref{LHsub}), and $G^{
%\ast
%}\in\RV(1-\gamma)$ with $\gamma\in\lbrack0,1),$ (\ref{HS})
%holding with $%
%\alpha=\beta+\gamma-1>0$ is \textbf{equivalent} to $\pibar^+_X$ $%
%\in\RV(-\beta)$ or $\Pi_X \in\MDA(\Lambda),$ but we omit the proof
%of the
%converse, as it follows from our later results. But note that
%Proposition~\ref{eq} implies that, when $\pibar^+_X$ $\in\MDA(\Lambda),$ we could
%also take $a(u)=\int_u^\infty\pibar^+_X(y)\rmd y/\pibar^+_X(u)$.
%\end{rem}

\begin{pf}
Assume (\ref{L-inf}),
%% and (\ref{LHsub}),
and that \eqref{Gcon} holds with $\gamma\in{}[0,1)$.

%%So by \eqref{tv22}, $A_{H^*}\in\RV(\gamma)$.
The starting point is Vigon's \'{e}quation amicale, \eqref{Vig2},
which we
write as $\overline{\Pi}^+_X(u)=I(u)+\mathrm{d}_{H^*}n(u), $ with
%$I(u)$ =$\int_{(u,\infty)}\pibar_{H^*}(y+u)\Pi_{\mathcal{H}} (
%\mathrm{d}y).$
%Integration by parts gives
%e4.11 #&#
\begin{eqnarray}
\label{I+I} I(u) &=&\int_{(0,\infty)}\Pi_{\mathcal{H}} (u+
\mathrm{d}y)\int_{(y,\infty)}\Pi_{H^{\ast}}(\mathrm{d}z)=\int
_{(0,\infty)}\Pi _{H^{\ast
}}(\mathrm{d}z)\int
_{(0,z)}\Pi_{\mathcal{H}} (u+\mathrm{d}y)
\nonumber
\\
&=&\int_{(0,\infty)}\Pi_{H^{\ast}}\bigl(a(u)\,\mathrm{d}z\bigr)
\Pi_{\mathcal{H}} \bigl\{\bigl(u,u+a(u)z\bigr]\bigr\}
\nonumber
\\[-8pt]
\\[-8pt]
\nonumber
&=& \biggl( \int_{(0,K]}+\int_{(K,\infty)}
\biggr) \Pi_{H^{\ast
}}\bigl(a(u)\,\mathrm{d}z\bigr)\Pi_{\mathcal{H}} \bigl
\{\bigl(u,u+a(u)z\bigr]\bigr\}
\\
&=:& I_{1}(u)+I_{2}(u), \qquad\mathrm{say},\nonumber
\end{eqnarray}
where $K>0$. Recall the definition of $A_{H^*}$ in \eqref{Adef2}, and note
that
\[
u\overline{\Pi}_{H^*}(u)\leq\int_0^u
\overline{\Pi }_{H^*}(y)\,\mathrm{d} y=A_{H^*}(u),\qquad u>0,
\]
so we have by the regular variation of $A_{H^*}$
\[
\frac{ a(u)I_{2}(u)}{A_{H^*}(a(u))\overline{\Pi}_{\mathcal{H}}(u)} \leq \frac{a(u)\overline{\Pi}_{H^*}(Ka(u))}{A_{H^*}(a(u))} \leq\frac
{A_{H^*}(Ka(u))}{KA_{H^*}(a(u))} \sim
\frac{1}{K^{1-\gamma}}.
\]
Since $0\le\gamma<1$ it follows that
%
%e4.12 #&#
\begin{equation}
\lim_{K\rightarrow\infty}\limsup_{u\rightarrow\infty}
\frac
{a(u)I_{2}(u)}{A_{H^*}(a(u)) \overline{\Pi}_{\mathcal{H}}(u)}=0. \label{tv9}
\end{equation}
Now assume \eqref{HS}, in which we set $\alpha=\beta+\gamma-1>0$. By
\eqref{tv5} with $\overline{F}$ replaced by $\overline{\Pi}_\cH$,
this implies
%
%e4.13 #&#
\begin{equation}
\label{pdef} \frac{\Pi_\cH\{(u,u+a(u)z]\}}{\overline{\Pi}_\cH(u)}\rightarrow \int_{0}^{z}p(y)
\,\mathrm{d}y
\end{equation}
uniformly for $z\in{}[0,K]$, where $p(\cdot)$ is the limiting
density associated with $\Pi_\cH$, that is, $\operatorname{Par}(\beta+\gamma-1)$ in
case \textup{(i)}, or $\operatorname{Exp}(1)$ in case \textup{(ii)}. So the component $I_1(u)$ in \eqref
{I+I} satisfies
%
%e4.14 #&#
\begin{eqnarray}
\label{cI} I_{1}(u) &\sim&\overline{\Pi}_{\mathcal{H}}(u)\int
_{0}^{K}\Pi _{H^{\ast
}}\bigl(a(u)\,\mathrm{d}z
\bigr)\int_{0}^{z}p(y)\,\mathrm{d}y
\nonumber
\\
&=&\overline{\Pi}_{\mathcal{H}}(u)\int_{0}^{K}p(y)
\,\mathrm{d}y\int_{y}^{K}\Pi_{H^{\ast}}
\bigl(a(u)\,\mathrm{d}z\bigr)
\\
&=&\overline{\Pi}_{\mathcal{H}}(u)\int_{0}^{K}p(y)
\overline{\Pi }_{H^*}\bigl(a(u)y\bigr)\,\mathrm{d}y-\overline{
\Pi}_{\mathcal{H}}(u)\overline {\Pi }_{H^*}\bigl(a(u)K\bigr)\int
_{0}^{K}p(y)\,\mathrm{d}y.\nonumber
\end{eqnarray}
\textup{(a)} When $\gamma\in(0,1)$, $A_{H^{\ast}}\in\RV(\gamma)$ is
equivalent, by the
monotone density theorem (Bingham, Goldie and Teugels \cite{BGT},
Theorem~1.7.2, page~39), to $\overline{\Pi}_{\mathcal{H}^{\ast}}\in
\RV
(\gamma-1)$, and then $\overline{\Pi}_{\mathcal{H}^{\ast}}(x)\sim
\gamma x^{-1}A_{H^{\ast}}(x)$. So
%
%e4.15 #&#
\begin{equation}
\int_{0}^{K}p(y)\overline{\Pi}_{H^*}
\bigl(a(u)y\bigr)\,\mathrm{d}y \sim\frac
{\gamma
A_{H^*}(a(u))}{a(u)}\int_{0}^{K}p(y)y^{\gamma-1}
\,\mathrm{d}y,
\end{equation}
and by taking $u\to\infty$ then $K\to\infty$ in \eqref{cI} we conclude
%
%e4.16 #&#
\begin{equation}
\label{tv10} \lim_{K\rightarrow\infty}\lim_{u\rightarrow\infty}
\frac
{a(u)I_{1}(u)}{A_{H^*}(a(u))\overline{\Pi}_{\mathcal{H}}(u)}=\gamma \int_{0}^{\infty
}p(y)y^{\gamma-1}
\,\mathrm{d}y= \gamma E\bigl(C^{\gamma-1}\bigr).
\end{equation}
\textup{(b)} When $\gamma=0$, so that $A_{H^*}$ is slowly varying, we use the
feature that $\lim_{x\downarrow0}p(x)=p(0)>0$ in either case,
$\operatorname{Par}(\beta
-1+\gamma)$ or $\operatorname{Exp}(1)$, to argue, given arbitrary $\varepsilon>0$, the
existence of a $\delta_{\varepsilon}>0$ such that for all large
enough $u$
\begin{eqnarray*}
a(u)\int_{0}^{\delta_{\varepsilon}}p(y)\overline{
\Pi}_{\mathcal
{H}}\bigl(a(u)y\bigr)\,\mathrm{d}y &\leq& p(0) (1+
\varepsilon)A_{H^*}\bigl(\delta_{\varepsilon}a(u)\bigr) \sim p(0) (1+
\varepsilon)A_{H^*}\bigl(a(u)\bigr)
\end{eqnarray*}
and
\begin{eqnarray*}
a(u)\int_{0}^{\delta_{\varepsilon}}p(y)\overline{\Pi
}_{H^*}\bigl(a(u)y\bigr)\,\mathrm{d%
}y & \geq& p(0) (1-
\varepsilon)A_{H^*}\bigl(\delta_\varepsilon a(u)\bigr) \sim p(0) (1-
\varepsilon)A_{H^*}\bigl(a(u)\bigr).
\end{eqnarray*}
$A_{H^*}$ slowly varying implies $x\overline{\Pi}_{H^*}(x)=\mathrm{o}
(A_{H^*}(x) )$ as $x\to\infty$, so with $\delta_{\varepsilon}$
fixed we can
argue
%%use the uniform convergence of $A_{H^*}$ on $[\delta_{\varepsilon
%},K]$ to get%
\begin{eqnarray*}
\int_{\delta_{\varepsilon}}^{K}p(y)\overline{\Pi }_{H^*}
\bigl(a(u)y\bigr)\,\mathrm{d}y &=& \mathrm{o} \biggl(\frac{1}{a(u)}\int
_{\delta_\varepsilon}^K p(y)A_{H^*}\bigl(a(u)y\bigr)
\frac{\,\mathrm{d} y}{y} \biggr)
\\
&=& \mathrm{o} \biggl(\frac{A_{H^*}(a(u))}{a(u)} \biggr),
\end{eqnarray*}
and we deduce for $\gamma=0$ that
%
%e4.17 #&#
\begin{equation}
\label{tv11} \lim_{K\rightarrow\infty}\lim_{u\rightarrow\infty}
\frac
{a(u)I_{1}(u)}{A_{H^*}(a(u))\overline{\Pi}_{\mathcal{H}}(u)}=p(0).
\end{equation}
Thus, in all cases we have
%
%e4.18 #&#
\begin{equation}
I(u)\sim\frac{c(\gamma,\beta) A_{H^*}(a(u))\overline{\Pi
}_{\mathcal
{H}} ( u ) }{a(u)} \label{tv21}
\end{equation}
for a constant $c(\gamma,\beta)\in(0,\infty)$, which we can
evaluate as
follows.

%%= \gamma E(C^{\gamma-1})$ when $\gamma\in(0,1)$ and
%$c(\gamma,\beta)=c(0,\beta)=p(0)=\beta$
\textup{(a)} When $\gamma\in(0,1)$, in case \textup{(i)}
%
%e4.19 #&#
\begin{equation}
\label{cgbdef} c(\gamma,\beta)= \gamma E\bigl(C^{\gamma-1}\bigr)=\gamma(
\beta+\gamma -1)\int_0^\infty\frac{x^{{\gamma}-1}\,\mathrm{d} x}{(1+x)^{\beta+\gamma}}=
\frac
{\Gamma(\gamma+1)\Gamma(\beta)}{\Gamma(\beta+\gamma-1)}.
\end{equation}
[Note that the density $p(\cdot)$ here is the one associated with
$\mathcal{H}$, not $X^+$, that is, it is Pareto with parameter $\alpha
=\beta+\gamma-1$; see \eqref{pdef}.]

In case \textup{(ii)}
%
%e4.20 #&#
\begin{equation}
\label{cgbdef2} c(\gamma,\beta)= \gamma E\bigl(C^{\gamma-1}\bigr)=\Gamma({
\gamma}+1).
\end{equation}
\textup{(b)} When $\gamma=0$, $p(0)=1-\beta$ in case \textup{(i)}, and in case \textup{(ii)},
$p(0)=1$, so we set $c(0,\beta)=1-\beta$ in case \textup{(i)}, and $c(0,\beta
)=1$ in case \textup{(ii)}.

Now %, still assuming \eqref{HS} with $\alpha=\beta+\gamma-1$,
%%it remains to treat $n(u)$ when $\rmd_{H^*}>0.$
integrate (\ref{Vig2}) and use the estimate (\ref{tv21}) to get
%
%e4.21 #&#
\begin{eqnarray}
\label{intX} \int_{u}^{\infty}\overline{
\Pi}_X^+(y)\,\mathrm{d}y &=& \int_{u}^{\infty}I(v)\,\mathrm{d}v+
\mathrm{d}_{H^*}\overline{\Pi}_{\mathcal{H}}(u)
\nonumber
\\[-8pt]
\\[-8pt]
\nonumber
&\sim& c(\gamma,\beta) \int_u^\infty
\frac{A_{H^*}(a(v))\overline
{\Pi
}_\cH(v)}{a(v)} \,\mathrm{d} v +\mathrm{d}_{H^*}\overline{
\Pi}_\cH(u).
\end{eqnarray}
Assume in addition that $\overline{\Pi}_{\mathcal{H}}\in\RV
(1-\gamma
-\beta
) $. This together with $A_{H^*}\in\RV(\gamma)$ means that the product
$\overline{\Pi}_\cH A_{H^*} \in\RV(1-\beta)$. Then, taking $a(u)=u$
in this
case, \eqref{intX} gives %$\int_{u}^{\infty}v^{-1}
%\pibar_{\mathcal{H}}\left( v\right) \sim(1-\gamma-\beta)^{-1}
%\overline{\Pi_{\mathcal{H}} }(u)$.
%e4.22 #&#
\begin{equation}
\label{i1} \frac{1}{\overline{\Pi}_\cH(u) A_{H^*}(u)} \int_u^\infty
\overline {\Pi}%
_X^+(y)\,\mathrm{d}y \sim c(\gamma,\beta)\int
_1^\infty v^{-\beta
}\,\mathrm {d} v +
\frac{\mathrm{d}_{H^*}}{A_{H^*}(u)}.
\end{equation}
In either case, $A_{H^*}(\infty)=\infty$ or $A_{H^*}(\infty)<\infty$,
we can
use the monotone density theorem again to deduce from this that
$\overline{\Pi}_X^+\in\RV(-\beta)$, and hence that \eqref{pX} holds
with $a(u)=u$.\vspace*{1pt}

%%, further,
%%\ben
%%\frac{u\pibar_X^+(u)}{\pibar_\cH(u) A_{H^*}(u)}
%%\sim
%%c_\gamma\int_1^\infty v^{-\beta}\rmd v +\frac{\rmd_{H^*}}{A_{H^*}(
%\infty)},
%%\een
%%which gives (\ref{est}) with $c_{\gamma,\beta}$ as the RHS in this
%case.
%%%Alternatively, when $EH_{1}^{\ast}=A_{H^*}(\infty)<\infty$, so
%%%$EX_{1}\in(-\infty,0)$, then we have $G^{\ast}(u)\sim u/\nu^{
%\ast},$ both terms %%have the same asymptotic behaviour, and we again
%get (\ref{est}).
Alternatively, suppose $\overline{\Pi}_{\mathcal{H}}\in\operatorname{MDA}(\Lambda)$. In this case, \eqref{intX} gives
\begin{eqnarray*}
\int_{u+xa(u)}^{\infty}\overline{\Pi}_X^+(y)
\,\mathrm{d}y &\sim& c(\gamma,\beta) \int_{u+xa(u)}^\infty
\frac{A_{H^*}(a(v))\overline
{\Pi
}_\cH(v)}{a(v)} \,\mathrm{d} v +\mathrm{d}_{H^*}\overline{
\Pi}_\cH \bigl(u+xa(u)\bigr),\qquad x\ge0.
\end{eqnarray*}
Change variable by setting $v=u+v^{\prime}a(u)$ on the RHS. Since
$a(\cdot)$ is
self-neglecting, we have $a(v)=a(u+v^{\prime}a(u))\sim a(u)$, so by the
regular variation of $A_{H^*}$,
\[
\frac{A_{H^*}(a(v))}{a(v)} \sim\frac{A_{H^*}(a(u))}{a(u)},
\]
and since $\overline{\Pi}_{\mathcal{H}}\in\operatorname{MDA}(\Lambda)$,
\[
\overline{\Pi}_\cH(v)= \overline{\Pi}_\cH
\bigl(u+v^{\prime}a(u)\bigr)\sim \mathrm{e}^{-v^{\prime}} \overline{
\Pi}_\cH(u).
\]
Thus, for $x \ge0$
%
%e4.23 #&#
\begin{eqnarray}
\label{eX}&& \frac{1}{\overline{\Pi}_\cH(u)}\int_{u+xa(u)}^{\infty}
\overline {\Pi }_X^+(y)\,\mathrm{d}y\nonumber\\
&&\quad\sim c(\gamma,\beta) a(u)\int
_{x}^\infty\frac{A_{H^*}(a(v'))\overline
{\Pi
}_\cH(v')}{a(v')\overline{\Pi}_\cH(u)} \,\mathrm{d}
v^{\prime
}+\mathrm {d}_{H^*}\frac{\overline{\Pi}_\cH(u+xa(u))}{\overline{\Pi}_\cH(u)}
\\
&&\quad \sim c(\gamma,\beta) A_{H^*}\bigl(a(u)\bigr)\int_{x}^\infty
\mathrm{e}^{-v^{\prime}} \,\mathrm{d} v^{\prime}+\mathrm{e}^{-x}
\,\mathrm{d}_{H^*},\nonumber
\end{eqnarray}
which, applied with $x=0$, also gives
\[
\frac{\int_{u+xa(u)}^{\infty}\overline{\Pi}_X^+(y)\,\mathrm
{d}y}{\int_{u}^{\infty}\overline{\Pi}_X^+(y)\,\mathrm{d}y} \to \mathrm{e}^{-x},\qquad  x\ge0.
\]
Applying Theorem~1.2.2(3) of de Haan and Ferreira \cite{HaanFerr}, we get
%Applying Thm 2.7.3(b) p.~110 of \cite{Haan} check! de Haan (1970) we
%get
\[
\frac{\overline{\Pi}_X^+(u+xa(u))}{\overline{\Pi}_X^+(u)} \to \mathrm{e}^{-x},\qquad x\ge 0,
\]
which is \eqref{pX} in this case, and this implies
%
%e4.24 #&#
\begin{equation}
\label{fX} \frac{\int_{u+xa(u)}^{\infty}\overline{\Pi}_X^+(y)\,\mathrm
{d}y}{a(u)\overline{\Pi}_X^+(u)} \to \mathrm{e}^{-x},\qquad x\ge0,
\end{equation}
hence
%
%e4.25 #&#
\begin{equation}
\label{bX} a(u)\sim\frac{\int_{u}^{\infty}\overline{\Pi}_X^+(y)\,\mathrm
{d}y}{\overline{\Pi}_X^+(u)},
\end{equation}
as claimed for this case.

It remains to prove \eqref{est}. In case \textup{(i)}, when $\overline{\Pi
}_\cH
\in
\RV(1-\gamma-\beta)$ and $\overline{\Pi}_X\in\RV(-\beta)$, the
relation \eqref{i1} gives
%
%e4.26 #&#
\begin{eqnarray}
\label{c1} \overline{\Pi}_X^+(u)&\sim&\frac{\beta-1}{u} \int
_u^\infty \overline{\Pi }_X^+(y)
\,\mathrm{d}y
\nonumber
\\[-8pt]
\\[-8pt]
\nonumber
&\sim& \biggl(c(\gamma,\beta)+\frac{(\beta
-1)\mathrm
{d}_{H^*}}{A_{H^*}(u)} \biggr)
\frac{\overline{\Pi}_{\mathcal
{H}}(u)A_{H^*}(u)}{u}.
\end{eqnarray}
\textup{(a)} When $\gamma\in(0,1)$, this implies \eqref{est} with $c_{\gamma
,\beta}=c(\gamma,\beta)+(\beta-1)\,\mathrm{d}_{H^*}/EH_1^*$, for
$EH_1^*\le\infty$. \textup{(b)} When $\gamma=0$, $c_{0,\beta}=c(0,\beta)$ for
$EH_1^*=\infty$ and, for $EH_1^*<\infty$.
\begin{eqnarray*}
c_{0,\beta} &=& c(0,\beta)+\frac{(\beta-1)\mathrm{d}_{H^\ast}}{EH_1^\ast
-\mathrm{d}_{H^\ast}} =\beta+\frac{(\beta-1)\mathrm{d}_{H^\ast}}{EH_1^\ast
-\mathrm{d}_{H^\ast}}
\\
&=&\frac{\beta EH_1^\ast-\mathrm{d}_{H^\ast}}{EH_1^\ast-\mathrm{d}_{H^\ast}} =\frac
{\beta
EH_{1}^\ast-\mathrm{d}_{H^\ast}}{A_{H^{\ast}}(\infty)}.
\end{eqnarray*}
In case \textup{(ii)}, when $\Pi_X\in\operatorname{MDA}(\Lambda)$, \eqref{eX} and
\eqref{fX} give, instead of \eqref{c1},
%
%e4.27 #&#
\begin{equation}
\label{c2} \overline{\Pi}_X^+(u)\sim\frac{1}{a(u)} \int
_u^\infty\overline {\Pi }_X^+(y)
\,\mathrm{d}y \sim \biggl(c(\gamma,\beta)+\frac{\mathrm
{d}_{H^*}}{A_{H^*}(a(u))} \biggr)
\frac{\overline{\Pi}_{\mathcal
{H}}(u)A_{H^*}(a(u))}{a(u)}.
\end{equation}
\textup{(a)} When $\gamma\in(0,1)$ this implies \eqref{est} with $c_{\gamma
,\beta}=c(\gamma,\beta)+\,\mathrm{d}_{H^*}/EH_1^*$, for $EH_1^*\le
\infty$.
\textup{(b)}~When $\gamma=0$, $c_{0,\beta}=1$ for $EH_1^*=\infty$ and,
for $EH_1^*<\infty$,
\begin{eqnarray*}
c_{0,\beta} &=&c(0,\beta)+\frac{\mathrm{d}_{H^{\ast}}}{EH_{1}^{\ast
}-\mathrm{d}_{H^{\ast}}} =1+\frac{\mathrm{d}_{H^\ast}}{EH_1^{\ast}-\mathrm{d}_{H^\ast}}
\\
&=&\frac{EH_1^\ast}{EH_1^\ast-\mathrm{d}_{H^\ast}} =\frac{EH_{1}^{\ast
}}{A_{H^{\ast
}}(\infty)}.
\end{eqnarray*}
This completes the proof of Proposition~\ref{eq}.
\end{pf}
Doney \cite{doneystf}, Corollary~4, page~31 (interchange $+/-$ in his
result), shows that,
when
$\lim_{t\to\infty}X_t= -\infty$ a.s., $E|X_{1}|<\infty$ if and
only if
$EH_1^*<\infty$, and then
$E|X_{1}|=qEH_1^*$. The following proposition
generalises this, allowing for $EH_1^*=\infty$.

%pr4.3 #&#
\begin{proposition}
\label{Q} Assume $\lim_{t\to\infty}X_{t}=-\infty$ a.s. and $A^{\ast
}_X(\infty)=\infty$, or, equivalently, $EH_1^*=\infty$.
%%equivalently, $A_{H^*}(\infty)=\infty$.
%%%, and that $A_{H^*}(x)\in\RV(\gamma)$, $\gamma>0$.
Then
%
%e4.28 #&#
\begin{equation}
\lim_{x\rightarrow\infty}\frac{A^{\ast}_X(x)}{A_{H^*}(x)}=q. \label{q}
\end{equation}
\end{proposition}

\begin{pf}
Assume $\lim_{t\rightarrow
\infty}X_{t}=-\infty$ a.s. and $A_{H^{\ast}}(\infty)=\infty$. The
integral term in (\ref{Vig3}) can be written as
\begin{eqnarray*}
&&\int_{(0,\infty)} \bigl( \overline{\Pi}_{H^{\ast}}(u)-\overline
{\Pi }_{H^{\ast}}(y+u) \bigr) \Pi_{\mathcal{H}}(\mathrm{d}y)\\
&&\qquad=\int
_{(0,\infty)}\overline{\Pi}_{\mathcal{H}}(y)\,\mathrm{d}_{y}
\bigl( \overline{\Pi}_{H^{\ast}}(u)-\overline{\Pi}_{H^{\ast}}(y+u) \bigr)
\end{eqnarray*}
after integrating by parts. So, by integrating (\ref{Vig3}) over
$1\leq
u\leq x$, we have
%
%e4.29 #&#
\begin{equation}
A_{X}^{\ast}(x)-q\int_{1}^{x}
\overline{\Pi}_{H^{\ast}}(u)\,\mathrm {d}u=\mathrm{d}_{\mathcal{H}}\bigl(
\overline{\Pi}_{H^{\ast
}}(1)-\overline {\Pi}_{H^{\ast}}(x)\bigr)+I(x),
\label{split}
\end{equation}
where
\[
I(x)=\int_{(0,\infty)}\Pi_{\mathcal{H}}(\mathrm{d}y)\int
_{1}^{x} \bigl( \overline{\Pi}_{H^{\ast}}(u)-
\overline{\Pi}_{H^{\ast}}(y+u) \bigr) \,\mathrm{d}u.
\]
We can bound the inner integral by
\begin{eqnarray*}
\biggl( \int_{1}^{x}-\int_{1+y}^{x+y}
\biggr) \overline{\Pi }_{H^{\ast
}}(u)\,\mathrm{d}u &=& \biggl( \int
_{1}^{1+y}-\int_{x}^{x+y}
\biggr) \overline{\Pi}_{H^{\ast}}(u)\,\mathrm{d}u
\\
&\leq&\int_{1}^{1+y}\overline{
\Pi}_{H^{\ast}}(u)\,\mathrm{d}u\leq y\overline{\Pi}_{H^{\ast}}(1).
\end{eqnarray*}
Then, for any $K>0$,
%
%e4.30 #&#
\begin{eqnarray}\label{Iu}
I(x) &\leq&\overline{\Pi}_{H^{\ast}}(1)\int_{0}^{K}y
\Pi_{\mathcal
{H}}(\mathrm{d}y)+\int_{K}^{\infty}
\Pi_{\mathcal{H}}(\mathrm {d}y)\int_{1}^{x}
\bigl( \overline{\Pi}_{H^{\ast}}(u)-\overline{\Pi }_{H^{\ast
}}(y+u)
\bigr) \,\mathrm{d}u
\nonumber
\\
&\leq&\overline{\Pi}_{H^{\ast}}(1)\int_{0}^{K}y
\Pi_{\mathcal
{H}}(\mathrm{d}y)+\int_{K}^{\infty}
\Pi_{\mathcal{H}}(\mathrm {d}y)\int_{1}^{x}
\overline{\Pi}_{H^{\ast}}(u)\,\mathrm{d}u
\\
&\leq&\overline{\Pi}_{H^{\ast}}(1)\int_{0}^{K}y
\Pi_{\mathcal
{H}}(\mathrm{d}y)+\overline{\Pi}_{\mathcal{H}}(K)A_{H^{\ast}}(x).
\nonumber
\end{eqnarray}
Since $A_{H^{\ast}}(\infty)=\infty$, when we divide by $A_{H^{\ast}}(x)$
and let $x\rightarrow\infty$ and then $K\rightarrow\infty$ in
\eqref{Iu},
we get $\lim_{x\rightarrow\infty}I(x)/A_{H^{\ast}}(x)=0$. Then
(\ref{q})
follows from \eqref{split}.
\end{pf}
%
% holds in this case.
% When $A_{H^*}(\infty)<\infty$, \eqref{Iu} shows that $\lim_{x\to
%\infty}I(x)=:I(\infty)<\infty$, so from \eqref{split},
% \[
% \lim_{x\to\infty}A^{\ast}_X(x)=q\int_{1}^{\infty}\pibar_{H^*}(u)\rmd
%u+\rmd_\cH\pibar_{H^*}(1)+I(\infty),
%\]
%so again (\ref{q}) holds.
%\halmos

%re4.1 #&#
\begin{rem}
\textup{(i)} We mention that a random walk version of Proposition~\ref{Q} is
(in a different notation) in Lemma~1 of Denisov, Foss, and Korshunov
\cite{DFK}.

\textup{(ii)} When \eqref{Gcon} holds, that is, $A_{H^*}\in\RV(\gamma)$ with
$\gamma\in[0,1)$, and $A_{H^*}(\infty)=\infty$, then $A_X^*(\infty
)=\infty$
and, by \eqref{q}, $A_X^*\in\RV(\gamma)$. The latter is equivalent to
%
%e4.31 #&#
\begin{equation}
\label{xpA} \lim_{x\to\infty}\frac{x\overline{\Pi}_X^-(x)}{A_X^*(x)}=\gamma.
\end{equation}
This is also true when $A_{H^*}(\infty)<\infty$, equivalently,
$A_X^*(\infty)<\infty$. [Compare with \eqref{xhA}.]
\end{rem}

%s5 #&#
\section{The case \texorpdfstring{$\gamma=0$}{gamma=0} (including finite mean)}\label{s4}

%%\label{secgameq0} %\begin{rem}
%Essentially the same proof applies in the case $\gamma=0$ and $\Pi_{
%\mathcal{H}} =\infty
%, $ and gives a "local" version of Thm 4.2.
%\end{rem}
Assume (\ref{L-inf}) and \eqref{Gcon} with $\gamma=0$, so
%% so we also have \eqref{tv22} with $\gamma=0$.
$A_{H^{\ast}}\in\RV(0)$, or, equivalently, $x\overline{\Pi
}_{H^{\ast
}}(x)=\mathrm{o}(A_{H^{\ast}}(x))$ as $x\rightarrow\infty$.
%%hence, by Proposition~\ref{Q}, $A_X^*(x)$, are slowly varying as $x\to
%\infty$.
Now (e.g., use Theorem~4.4 of Doney and Maller \cite{DM2002} with $+/-$
interchanged) (\ref{L-inf}) implies
%
%e5.1 #&#
\begin{equation}
\frac{x\overline{\Pi}^{+}(x)}{A_{X}^{\ast}(x)}\leq\frac
{A_{X}^{+}(x)}{A_{X}^{\ast}(x)}\rightarrow0 \qquad\mbox{as } x\rightarrow
\infty\label{Apm}
\end{equation}
if $A_{X}^{\ast}(\infty)=\infty$, otherwise $A_{X}^{\ast}(\infty
)<\infty$ and then $A_{X}^{+}(\infty)<\infty$ and $\lim_{x\rightarrow
\infty}x\overline{\Pi}^{+}(x)=0$. Thus, since also $x\overline{\Pi
}_{X}^{-}(x)=\mathrm{o}(A_{X}^{\ast}(x))$ by \eqref{xpA},
\[
A(x):=\gamma+\overline{\Pi}^{+}(1)-\overline{\Pi }^{-}(1)+A_{X}^{+}(x)-A_{X}^{\ast}(x)
\sim-A_{X}^{\ast}(x)\qquad \mbox{as } x\rightarrow\infty,
\]
and we see that $x\overline{\Pi}(x)=\mathrm{o}(-A(x))$ as $x\rightarrow\infty$.
This means that $X_{t}$ is negatively relatively stable (Doney and
Maller \cite{DM2002}), or, equivalently,
%%$A(x)$ is slowly varying as $x\to\infty$.%%, Hence Thus
$X_{t}^{\ast}$ is positively relatively stable, as $t\rightarrow
\infty$.
Consequently, we can employ a version of the weak law of large numbers even
if the mean is infinite; specifically there is a continuous, increasing
function $c(\cdot)\in\RV(1)$ such that $X_{t}^{\ast}/c(t)\overset
{\mathrm{P}}{\longrightarrow}1$ as $t\rightarrow\infty$. The
function $c(\cdot)$
can be chosen to be strictly increasing and to satisfy
\[
c(x)=xA_{X}^{\ast}\bigl(c(x)\bigr),\qquad x>0,
\]
and its inverse function $b(\cdot):=c^{-1}(\cdot)$ is given by
\[
b(y)=\frac{y}{A_{X}^{\ast}(y)},\qquad y>0.
\]
Employing Proposition~\ref{Q}, we see that
%
%e5.2 #&#
\begin{equation}
b(y)=\frac{y}{A_{X}^{\ast}(y)}\sim\frac{y}{qA_{H^{\ast}}(y)}\qquad\mbox{as } y\rightarrow\infty,
\label{a2}
\end{equation}
when $A_{H^{\ast}}(\infty)=\infty$. When $A_{H^{\ast}}(\infty
)<\infty$, and so $EX_{1}\in(-\infty,0)$, we simply take $c(x)=|EX_{1}|x$
and $b(x)=x/|EX_{1}|$, $x>0$.

We define another norming function by $r(u)=b(a(u))$, and note that
$c(r(u))=a(u)$ and
%
%e5.3 #&#
\begin{equation}
\label{a9} r(u)\sim\frac{a(u)}{qA_{H^*}(a(u))}
\end{equation}
when $A_{H^*}(\infty)=\infty$, and
%
%e5.4 #&#
\begin{equation}
\label{a9a} r(u)=\frac{a(u)}{|EX_1|}= \frac{a(u)}{qEH_1^*}
\end{equation}
when $A_{H^*}(\infty)<\infty$. The function $r(u)$ turns out to be the
right norming for $\tau_u$ in the present situation.

\begin{pf*}{Proof of Theorem \ref{old}}
Assume (\ref{L-inf})
and (\ref{LHsub}), and that \eqref{Gcon} holds with $\gamma=0$. Then
parts~1\textup{(a)} and \textup{(b)} of the theorem are equivalent by Proposition~\ref{nasc}
applied to the subordinator $\mathcal{Y}:=\mathcal{H}$, and part 1\textup{(c)}
follows from part~1\textup{(b)} by Proposition~\ref{eq}. We now show that part 1\textup{(c)}
implies part~2.
\end{pf*}

%pr5.1 #&#
\begin{proposition}\label{Knew}
Assume (\ref{L-inf}) and (\ref{LHsub}), and additionally that
$A_{H^*}\in\RV(0)$, and either \textup{(i)}~$\overline{\Pi}^+_X\in\RV
(-\beta
)$, where
$\beta>1$, or \textup{(ii)} $\Pi_X\in\operatorname{MDA}(\Lambda)$. Then the
conclusions of part 2 of Theorem~\ref{old} hold.
\end{proposition}

\begin{pf}
A slight extension of a result proved in Doney and Rivero \cite{DR2012}
states that, on the event $X_{\tau_u-}<u$,
the joint distribution of $(\tau_u,X_{\tau_u-})$ is given by
%
%e5.5 #&#
\begin{eqnarray}
\label{RADVR} &&P ( \tau_u\in\mathrm{d} t,X_{\tau_u-}\in
\mathrm{d} y )
\nonumber
\\[-8pt]
\\[-8pt]
\nonumber
&&\quad =P(X_t\in\mathrm{d} y,\overline{X}_t\leq
u)\overline{\Pi }^+_X(u-y)\,\mathrm{d} t,\qquad t>0, u>0, y\in\R.
\end{eqnarray}
Thus, $\tau_u$ has a density, and for $\veps>0$ we can write
(recall that $Z^{(u)}=-X_{\tau_u-}=X_{\tau_u-}^*$)
\begin{eqnarray*}
&&P\bigl(\tau_u \in r(u)\,\mathrm{d} t,Z^{(u)}\in{}\bigl[(1-
\veps)c\bigl(\tau (u)\bigr),(1+\veps)c(\tau_u)\bigr]\bigr)
\\
&&\quad = \int_{[(1-\veps)c({tr}(u)), (1+\veps)c(tr(u))]}\overline {\Pi }^+_X(u+y) P
\bigl(X_{tr(u)}^{\ast}\in\mathrm{d} y,\overline
{X}_{tr(u)}\leq u\bigr)\,\mathrm{d} t.
\end{eqnarray*}
Under the assumptions of the proposition, the limit relation \eqref{pX}
holds, and also $c(\cdot)\in\RV(1)$ implies $c(tr(u))\sim
t(c(r(u)))=ta(u)$. So the last
integral is asymptotically equivalent to
\begin{eqnarray*}
&&\int_{[(1-\veps)t,(1+\veps)t]}\overline{\Pi }^+_X\bigl(u+ya(u)
\bigr)P\bigl(X_{tr(u)}^{\ast}\in a(u)\,\mathrm{d} y,\overline
{X}_{tr(u)}\leq u\bigr)
\\
&&\quad\sim r(u)\overline{\Pi}^+_X(u)\int_{[(1-\veps)t,(1+\veps)t]}P(C>y)
P\bigl(X_{tr(u)}^{\ast}\in a(u)\,\mathrm{d} y,\overline{X}_{tr(u)}
\leq u\bigr)
\\
%&=& r(u)\overline{\Pi}^+_X(u)\int_{[(1-\veps)t,(1+\veps)t]}P(C>y)
%P\left(\frac{X_{tr(u)}^{\ast}}{c(r(u))}\in\,\mathrm{d} y,\frac{
%\overline{X}%
%_{tr(u)}}{c(r(u))}\leq\frac{u}{a(u)}\right) \\
&&\quad= r(u)\overline{
\Pi}^+_X(u) \biggl\{ \int_{[1-\veps,1+\veps]}P(C>ty) P \biggl(
\frac{X_{tr(u)}^{\ast}}{ta(u)}\in\mathrm{d} y \biggr)+\mathrm{o}(1) \biggr\},
\end{eqnarray*}
where we use the fact that $P(\overline{X}_{tr(u)}>u)\le P(\overline
{X}_\infty>u)\rightarrow0$ as $u\to\infty$.
This follows because $\overline{X}_\infty=\sup_{t\ge0}X_t$ is a finite
r.v. a.s. under \eqref{L-inf}.

Next, since
\[
\frac{X_{tr(u)}^{\ast}}{ta(u)} \sim \frac{X_{tr(u)}^{\ast}}{tc(r(u))}\topr1,
\]
for all $t>0$, we deduce that
\[
\int_{[1-\veps,1+\veps]}P(C>ty)P \bigl(X_{tr(u)}^*\in ta(u)
\,\mathrm{d} y \bigr)=P(C>t)+\mathrm{o}(1),
\]
so that
%
%e5.6 #&#
\begin{eqnarray}
\label{5.3a} &&P^{(u)} \bigl(\tau_u \in r(u)\,\mathrm{d}
t,Z^{(u)}\in \bigl[(1-\veps )c(\tau_u),(1+\veps)c(
\tau_u)\bigr] \bigr)
\nonumber
\\
& &\quad\sim \frac{r(u)\overline{\Pi}^+_X(u)P(C>t)\,\mathrm{d} t}{P(\tau
_u<\infty)}
\nonumber
\\[-8pt]
\\[-8pt]
\nonumber
& &\quad\sim \frac{a(u)\overline{\Pi}^+_X(u)P(C>t)\,\mathrm{d}
t}{\overline
{\Pi}_{\mathcal{H}}(u)A_{H^*}(u)} \qquad\bigl(\mbox{by \eqref{max} and \eqref{a9}}\bigr)
\nonumber
\\
&&\quad\to c_{0,\beta}P(C>t)\,\mathrm{d} t \qquad\bigl(\mbox{by \eqref{est}}\bigr).\nonumber
\end{eqnarray}
The evaluation of $c_{0,\beta}$ from \eqref{est}
(and see the end of the proof of Proposition~\ref{eq})
shows that the limit here is a probability
density function, and since it does not depend on $\varepsilon$, we deduce
that (\ref{local}) holds, and also that, conditioned on $\tau_u=tr(u)$,
the $P^{(u)}$-distribution of $X^{\ast}(\tau_u-)/c(\tau_u)$ converges
to the
distribution concentrated on $1$.
\end{pf}

%re5.1 #&#
\begin{rem}\label{ont}
The event $\{X_{\tau_u-}<u\}$ to which \eqref{RADVR} is restricted has
$P^{(u)}$-probability approaching 1 as $u\to\infty$.
This follows since $\lim_{u\to\infty}P(Z^{(u)}/a(u)\le0)=0$
in conditions \textup{(a)}--(c) of Theorem~\ref{old} (and similarly in Theorem~\ref{new}), so we have
\[
P^{(u)}(X_{\tau(u)-}=u)=P^{(u)}\bigl(Z^{(u)}=-u
\bigr) \leq P^{(u)}\bigl(Z^{(u)}\leq0\bigr)\rightarrow0 \qquad\mbox{as } u\to\infty.
\]
\end{rem}

To extend (\ref{local}) to the $k$-dimensional distributions, we take
$0<s_{1}<s_{2}<\cdots<s_{k}<1$, set
\[
A_{k}:= \biggl\{ 1-\varepsilon\leq\frac{X^{\ast}(s_{i}\tau
_u)}{s_{i}c(\tau_u)}\leq1+\varepsilon
\mbox{ for } i=1,2,\ldots ,k \biggr\},
\]
and apply the previous argument to
\[
P\bigl(A_{k},\tau_u\in r(u)\,\mathrm{d} t,Z^{(u)}
\in{}\bigl[(1-\veps )c\bigl(\tau (u)\bigr),(1+\veps)c(\tau_u)\bigr]
\bigr).
\]
We find that
\begin{eqnarray*}
P^{(u)}\bigl(A_{k},\tau_u \in r(u)\,\mathrm{d}
t,Z^{(u)}\in{}\bigl[(1-\veps )c(\tau_u),(1+\veps)c(
\tau_u)\bigr]\bigr)\rightarrow c_{0,\beta}P(C>t)\,\mathrm{d} t,
\end{eqnarray*}
and the convergence of the $k$-dimensional distributions follows.
%%%\halmos

To include the behaviour of the overshoot, we need the following result.

%le5.1 #&#
\begin{lemma}
\label{Lovun} For $u>0$, $z\ge0$, and $x\geq0$ we have
\[
P^{(u)}\bigl(Z^{(u)}\in\,\mathrm{d}z,\mathrm{O}^{(u)}>x
\bigr)=P^{(u)}\bigl(Z^{(u)}\in\mathrm {d}z\bigr)%
\frac{\overline{\Pi}^+_X(u+x+z)}{\overline{\Pi}^+_X(u+z)}.
\]
\end{lemma}

\begin{pf}
Using the quintuple law in Doney and Kyprianou \cite{DK} twice %%(see
%also Griffin and Maller (2011)),
gives
\begin{eqnarray*}
P\bigl(Z^{(u)}\in\mathrm{d}z,\mathrm{O}^{(u)}>x\bigr) &=& \int
_{0<w\le u}G(\mathrm{d}w)G^{\ast}(u-w-\mathrm{d} z)\overline{\Pi
}^+_X(u+x+z)
\\
&=& \int_{0<w\le u}G(\mathrm{d}w)G^{\ast}(u-w-\mathrm{d} z)
\overline{\Pi }^+_X(u+z)\frac{\overline{\Pi}_X^+(u+x+z)}{\overline{\Pi
}^+_X(u+z)}
\\
&=&P^{(u)}\bigl(Z^{(u)}\in\mathrm{d}z\bigr)\frac{\overline{\Pi
}^+_X(u+x+z)}{\overline{\Pi}_X^+(u+z)}.
\end{eqnarray*}
(Note that there is no issue of creeping to take into account since
$\mathrm{O}^{(u)}>0$ implies $X_{\tau_{u}}>u$.)
\end{pf}

%co5.1 #&#
\begin{corollary}
\label{Kcor} Under the assumptions of Proposition~\ref{Knew}, the
$P^{(u)}$-finite-dimensional distributions $\mathbf{Y}^{(u)}$, defined
in (\ref{Y}),
converge to those of $ ( V,U,V,(V\mathbf{D}^{(0)}(s))_{0\le s\le
1} )$.
\end{corollary}

\begin{pf}The result for
\[
\biggl( \frac{Z^{(u)}}{a(u)}, \frac{\tau_u}{b(a(u))}, \biggl( \frac
{X^{\ast
}(s\tau_u)}{a(u)}
\biggr) _{0\le s\le1} \biggr)
\]
is immediate from Proposition~\ref{Knew}, and since, given $Z^{(u)}$,
$\mathrm{O}^{(u)}$ is independent of the pre-$\tau_u$ $\sigma$-field, we need
only check that
\[
P\bigl(\mathrm{O}^{(u)}>xa(u)\mid Z^{(u)}=a(u)z\bigr)\rightarrow\cases{
\displaystyle\biggl( \frac{1+z}{1+z+x} \biggr)
^\beta, &\quad $\mbox {in case \textup{(i)},}$
\vspace*{2pt}\cr
\mathrm{e}^{-x}, & \quad$\mbox{in case \textup{(ii)}.}$}
\]
But this is immediate from Lemma~\ref{Lovun}.
\end{pf}

In particular, when part 1\textup{(c)} of Theorem~\ref{old} holds, we have from
Corollary~\ref{Kcor} that the $P^{(u)}$-distribution of $\mathrm{O}^{(u)}$ converges
to that of $U$, so 1\textup{(a)}
%, and hence the assumption ofProposition~\ref{nasc} holds, so 1\textup{(b)} also
holds. Thus, parts 1\textup{(a)}--\textup{(c)} are proved equivalent.

Finally, for part~3 of Theorem~\ref{old}, we show that the convergence in
this result can be replaced by weak convergence on the Skorokhod space.

%pr5.2 #&#
\begin{proposition}
\label{wh?} Under the assumptions of Proposition~\ref{Knew}, the
$P^{(u)}$-distribution of $\mathbf{Y}^{(u)}$ converges weakly on
$\mathbb{R}^{3}\times\mathbb{D}_{0}[0,1]$ as $u\rightarrow\infty$.
\end{proposition}

\begin{pf}
Put $\mathbf{Y}^{(u)}=(W^{(u)},\mathbf{X}^{(u)})$, where
\[
W^{(u)}:= \biggl( \frac{Z^{(u)}}{a(u)}, \frac{\mathrm{O}^{(u)}}{a(u)},
\frac
{\tau_u}{b(a(u))} \biggr)\quad \mathrm{and}\quad \mathbf {X}^{(u)}:= \biggl(
\frac{X^{\ast}(s\tau_u)}{a(u)} \biggr) _{0\le s\le1}.
\]
We need only prove tightness. This will follow if we can show that for
any $%
\varepsilon>0$ there is a compact subset of $K$ of $\mathbb
{R}^{3}\times
\mathbb{D}_0[0,1]$ such that $\lim\sup_{u\rightarrow\infty}P^{(u)}(
\mathbf{Y^{(u)}}\in K^{c})\leq\varepsilon$. We will do this with $%
K=K_{1}\times K_{2}$, where $K_{1}\subset\mathbb{R}^{3}$ is of the
form $%
\{1/D<x_{r}<D,r=1,2,3\}$, $K_{2}\subset\mathbb{D}_0[0,1]$ will be specified
later, and $D$ is fixed with $P^{(u)}(W^{(u)}\in K_{1}^{c})\leq
\varepsilon
/2$ for large $u$. So it suffices to show that $\lim\sup_{u\rightarrow
\infty}P^{(u)}(\mathbf{Y}^{(u)}\in K_{1}\times K_{2}^{c})\leq
\varepsilon
/2$. This probability is dominated by
\[
P^{(u)} \bigl( B\cap\bigl(\mathbf{X}^{(u)}\in
K_{2}^{c}\bigr) \bigr)
\]
where
\[
B= \biggl\{ \frac{\tau_u}{r(u)}\in\bigl(D^{-1},D\bigr),
\frac
{Z^{(u)}}{a(u)}\in \bigl(D^{-1},D\bigr) \biggr\}.
\]
But (recall $c(r(u))=a(u)$ and \eqref{RADVR})
%
%e5.7 #&#
\begin{eqnarray}
\label{qery} &&P^{(u)} \bigl(\mathbf{X}^{(u)} \in
K_{2}^{c},B \bigr)
\nonumber\\
&&\quad\le \frac{1}{P(\tau_u<\infty)} \int
_{r(u)/D}^{r(u)D}\int_{z\in
(D^{-1},D)}\,\mathrm{d}
tP \bigl( X_{t}^{\ast}\in a(u)\,\mathrm{d} z,\mathbf
{X}^{(u)}\in K_{2}^{c} \bigr) \overline{
\Pi}^+_X\bigl(u+a(u)z\bigr)
\nonumber
\\[-8pt]
\\[-8pt]
\nonumber
&&\quad\leq \frac{\overline{\Pi}^+_X(u)}{P(\tau_u<\infty)} \int
_{r(u)/D}^{r(u)D}\,\mathrm{d} tP \biggl( \biggl(
\frac{X_{st}^{\ast
}}{a(u)}, 0\le s\le1 \biggr)\in K_{2}^{c} \biggr)
\\
&&\quad= \frac{r(u)\overline{\Pi}^+_X(u)}{P(\tau_u<\infty)} \int_{1/D}^{D}\,\mathrm{d}
tP \biggl( \biggl( \frac{X_{r(u)st}^{\ast
}}{c(r(u))}, 0\le s\le1 \biggr) \in
K_{2}^{c} \biggr).\nonumber
\end{eqnarray}
As shown in (\ref{5.3a}), the factor
\[
\frac{r(u)\overline{\Pi}^+_X(u)}{P(\tau_u<\infty)}\to c_{0,\beta
}\qquad \mbox{as } u\to\infty.
\]
Also, since $(X_{ys}^{\ast}/c(y))_{0\le s\le1}$ is tight as
$y\rightarrow
\infty$, we can choose $K_{2} $ such that when $D^{-1}a(u)$ is sufficiently
large,
\[
P \biggl( \sup_{t\in(D^{-1},D)} \biggl( \frac{X_{r(u)st}^{\ast
}}{c(r(u))}, 0\le s\le1
\biggr) \in K_{2}^{c} \biggr) \leq \varepsilon,
\]
and the result follows.
\end{pf}

%s6 #&#
\section{The case \texorpdfstring{$0<\gamma<1$}{0<gamma<1} (infinite mean)}\label{s5}

Throughout this section, our standing assumptions (and notations) will
be those of Theorem~\ref{new}, namely, (\ref{L-inf}) and (\ref{LHsub})
hold, and \eqref{Gcon}
holds with $\gamma\in(0,1)$.
%Recall from (\ref{tv22}) that $G^{\ast}\in\RV(1-\gamma)$ is
%equivalent to $A_{H^*}(x)\sim k_\gamma x/G^*(x)\in\RV(\gamma),$ and
%hence, by Proposition~\ref{Q}, to
%$A^{\ast}_X(x)\sim qk_{\gamma}x/G^{\ast}(x)$ as $x\to\infty$.
By the monotone density theorem, the latter is equivalent to
%
%e6.1 #&#
\begin{equation}
\overline{\Pi}^-_X(x)\sim\gamma x^{-1}A^{\ast}_X(x)
\in\RV(\gamma -1) \qquad\mbox{as } x\rightarrow\infty. \label{Gstar}
\end{equation}
From \eqref{Apm}, we then deduce $\lim_{x\to\infty}\overline{\Pi
}^+_X(x)/\overline{\Pi}^-_X(x)=0$. This together with \eqref{Gstar}
means that $X^{\ast}$ is in the domain of attraction of a standard
stable subordinator,
$\mathbf{D}$, of parameter $\overline{\gamma}:=1-\gamma\in(0,1)$.
Thus, we
can find a continuous, increasing function $c(\cdot)$ such that $
(X_{su}^{\ast}/c(u) ) _{s>0}\overset{D}{\rightarrow}\mathbf{D,}$
and one can check that
\[
u\overline{\Pi}^-_X\bigl(c(u)\bigr)\rightarrow1/\Gamma(\gamma).
\]
Write $b(\cdot)$ for the inverse of $c(\cdot)$, so that $b(\cdot)\in
\mathcal{R}_{\overline{\gamma}}$, and
%
%e6.2 #&#
\begin{equation}
b(u)\sim\frac{1}{\Gamma(\gamma)\overline{\Pi}^-_X(u)}. \label{y}
\end{equation}
Put $r(u)=b(a(u))$, so that
%
%e6.3 #&#
\begin{equation}
\label{g} r(u)\sim\frac{1}{\Gamma(\gamma)\overline{\Pi}^-_X(a(u))} \sim \frac
{a(u)}{\Gamma(1+\gamma)A_X^*(a(u))} \qquad\bigl(\mbox{by }\eqref{Gstar}\bigr).
\end{equation}

A version of Stone's stable local limit theorem (see Proposition~13 of
Doney and Rivero \cite{DR2012}) implies that
%
%e6.4 #&#
\begin{equation}
P(X_{{tv}}^{\ast}\in\bigl(c(v)z,c(v)z+\Delta]\bigr)=
\frac{\Delta}{c(v)} \bigl( h_{t} ( z ) +\mathrm{o}(1) \bigr) \label{llt+}
\end{equation}
as $v\rightarrow\infty$, uniformly for $z\in\mathbb{R}$,
$\Delta\in{}[\Delta_{0},\Delta_{1}]$, for any fixed $0<\Delta
_{0}<\Delta_{1}<\infty$, and $t\in{}[ T_{0},T_{1}]$, for any
fixed $0<T_{0}<T_{1}<\infty$. Here $h_{t}(z)\,\mathrm{d}z=P(D_{t}\in
\mathrm{d}z)$
[see \eqref{hdef}], so that, in particular, the term $h_{t} (
z ) $
is zero for $z<0$. A simple consequence of this is the existence of
constants $v_{0}$ and $C$ such that for all $v\geq v_{0,}$ $\Delta\in
{}[\Delta_{0},\Delta_{1}]$, and $t\in{}[ T_{0},T_{1}]$,
%
%e6.5 #&#
\begin{equation}
P\bigl(X_{{tv}}^{\ast}\in\bigl(c(v)z,c(v)z+\Delta\bigr]\bigr)\leq
\frac{C\Delta}{c(v)}. \label{llb}
\end{equation}
Notice that if we put $v=r(u)$ in (\ref{llt+}) we have
$c(v)=c(b(a(u)))=a(u)$, so an equivalent version of (\ref{llt+}) is
%
%e6.6 #&#
\begin{equation}
P\bigl(X_{tr(u)}^{\ast}\in\bigl(a(u)z,a(u)z+\Delta\bigr]\bigr)=
\frac{\Delta
}{a(u)} \bigl( h_{t} ( z ) +\mathrm{o}(1) \bigr)\qquad \mbox{as }u
\rightarrow\infty. \label{llt++}
\end{equation}

We have already proved part 1 of Theorem~\ref{new}, except for the
implication from parts 1\textup{(c)} to~\textup{(a)}, and we now show that part 1\textup{(c)}
implies part 2, and then that this implies part 1\textup{(a)}.

%pr6.1 #&#
\begin{proposition}
\label{W} Assume (\ref{L-inf}) and (\ref{LHsub}), and that
$A_{H^*}\in
\RV(\gamma)$ with $\gamma\in(0,1)$. Suppose either \textup{\textup{(i)}} $\overline
{\Pi
}^+_X(x)\in\RV(-\beta)$, where $\beta>1-\gamma$, or \textup{\textup{(ii)}}
$\overline
{\Pi}^+_X(x)\in\operatorname{MDA}(\Lambda)$ and $\overline{\Pi
}_{\mathcal{H}}\in
\mathcal{S}$. Then part~2 of Theorem~\ref{new} holds.
\end{proposition}

\begin{pf}Under the conditions of
the proposition, we have from \eqref{RADVR}
%
%e6.7 #&#
\begin{eqnarray}
\label{6.6a} && P\bigl(\tau_{u}\in r(u)\,\mathrm{d}t,Z^{(u)}
\in{}\bigl[ za(u),za(u)+\Delta \bigr]\bigr)\nonumber
\\
&&\quad=\int_{y\in{}[0,\Delta]}
\overline{\Pi} _{X}^{+}\bigl(u+za(u)+y\bigr)P
\bigl(X_{tr(u)}^{\ast}\in za(u)+\mathrm{d}y,\overline
{X}_{tr(u)}\leq u\bigr)\,\mathrm{d}t
\nonumber
\\[-8pt]
\\[-8pt]
\nonumber
&&\quad\sim \overline{
\Pi}_{X}^{+}\bigl(u+za(u)\bigr)\int_{y\in{}[0,\Delta
]}P
\bigl(X_{tr(u)}^{\ast}\in za(u)+\mathrm{d}y,\overline{X}_{tr(u)}
\leq u\bigr)\,\mathrm{d}t
\\
&&\quad\sim \overline{\Pi}_{X}^{+}(u)P(C>z)P
\bigl(X_{tr(u)}^{\ast}\in{}\bigl[ za(u),za(u)+\Delta\bigr],
\overline{X}_{tr(u)}\leq u\bigr)\,\mathrm{d}t.\nonumber
\end{eqnarray}
Write
\[
P\bigl(X_{tr(u)}^{\ast}\in{}\bigl[ za(u),za(u)+\Delta\bigr],
\overline {X}_{tr(u)}\leq u\bigr)=P_{1}(u)-P_{2}(u),
\]
where, by (\ref{llt+}),
%
%e6.8 #&#
\begin{equation}
P_{1}(u):=P\bigl(X_{tr(u)}^{\ast}\in{}\bigl[
za(u),za(u)+\Delta\bigr]\bigr)=\frac
{\Delta}{a(u)} \bigl( h_{t} ( z )
+\mathrm{o}(1) \bigr), \label{e1}
\end{equation}
and we will show that
%
%e6.9 #&#
\begin{equation}
P_{2}(u):=P\bigl(X_{tr(u)}^{\ast}\in{}\bigl[
za(u),za(u)+\Delta\bigr],\overline {X}_{tr(u)}>u\bigr)=\mathrm{o} \biggl(
\frac{\Delta}{a(u)} \biggr),\qquad u\to\infty. \label{f}
\end{equation}
To do this, observe that $\{\overline{X}_{tr(u)}>u\}\subseteq
\{\tau_u\le tr(u)\}$, and decompose $P_2(u)$ further according as
$\tau_u\le tr(u)/2$ or $tr(u)/2<\tau_u\le tr(u)$.
Thus, write $P_{2}(u)=P_{2}^{(1)}(u)+P_{2}^{(2)}(u)$, recall that
$\mathrm{O}^{(u)}$ is independent of the pre-$\tau_u$ $\sigma$-field, and argue
as follows: %%%% where, writing $tr(u)=r,$
%e6.10 #&#
\begin{eqnarray}
\label{6.8a} P_{2}^{(1)}(u)&:=&P \bigl(\tau_{u}
\leq tr(u)/2,X_{tr(u)}^{\ast}\in {}\bigl[ za(u),za(u)+\Delta\bigr]
\bigr)
\nonumber
\\
&=&\int_{0\leq s\leq tr(u)/2}\int_{x>0}P\bigl(
\tau_{u}\in \mathrm{d}s,\mathrm{O}^{(u)}\in \,\mathrm{d}x\bigr)
\nonumber
\\
&&\hspace*{67pt}{} \times P\bigl(X_{tr(u)-s}^{\ast}\in{}\bigl[
u+x+za(u),u+x+za(u)+\Delta\bigr]\bigr)\qquad\qquad
\\
\nonumber&\leq&\int_{0\leq s\leq tr(u)/2}\int_{x>0}P\bigl(
\tau_{u}\in \mathrm{d}s,\mathrm{O}^{(u)}\in \mathrm{d}x\bigr)\frac{C\Delta}{c(tr(u)-s)}
\qquad\bigl(\mbox{by (\ref{llb})}\bigr)
\nonumber
\\
&\leq&\frac{C^{\prime}\Delta}{c(tr(u))}P(\tau_{u}<\infty)
\nonumber\\
&=&\mathrm{o} \biggl(
\frac{\Delta}{c(tr(u))} \biggr).\nonumber
\end{eqnarray}
Next, introduce $\tau^{\ast}(u)=\inf\{s\dvt X_{s}^{\ast}>u\}$ and
$\sigma
_{v}(u)=\sup\{s\leq v\dvt X_{s}>u\}$. Use the duality lemma (Bertoin \cite
{bert}, page~45)
to see that for any $w$ and any $v>0$
\[
P\bigl(\sigma_{v}(u)\in \mathrm{d}s\mid X_{v}^{\ast}=w
\bigr)=P\bigl(\tau^{\ast}(u+w)\in v-\mathrm{d}s\mid X_{v}^{\ast}=w
\bigr).
\]
Applying this with $v=tr(u)$ and $w=za(u)+y$ gives
\begin{eqnarray*}
P_{2}^{(2)}(u) &=&\int_{[0,\Delta]}P \bigl(
tr(u)/2<\tau_{u}\leq tr(u), X_{tr(u)}^{\ast}\in za(u)+
\mathrm{d}y \bigr)
\\
&\leq&\int_{[0,\Delta]}P \bigl( tr(u)/2<\sigma_{tr(u)}(u)
\leq tr(u), X_{tr(u)}^{\ast}\in za(u)+\mathrm{d}y \bigr)
\\
&=&\int_{[0,\Delta]}P \bigl( 0<\tau^{\ast
}\bigl(u+za(u)+y
\bigr)<tr(u)/2, X_{tr(u)}^{\ast}\in za(u)+\mathrm{d}y \bigr)
\\
&\leq&P \bigl( 0<\tau^{\ast}\bigl(u+za(u)\bigr)<tr(u)/2,
X_{tr(u)}^{\ast
}\in \bigl(za(u),za(u)+\Delta] \bigr)
\\
&=&\int_{0\leq v\leq tr(u)/2}\int_{y>0}P\bigl(
\tau^{\ast}\bigl(u+za(u)\bigr)\in \mathrm{d} v, X_{v}^{\ast}
\in u+za(u)+\mathrm{d}y\bigr)
\\
&&\hspace*{68pt}{}\times P\bigl(X_{tr(u)-v}\in(u+y-\Delta,u+y]\bigr)
\\
&=&\mathrm{o}(1)\int_{0\leq v\leq tr(u)/2}P\bigl(\tau^{\ast}\bigl(u+za(u)
\bigr)\in\mathrm{d}v\bigr) \frac{\Delta}{c(tr(u)-v)}
\\
&=&\mathrm{o} \biggl( \frac{\Delta}{c(tr(u))} \biggr).
\end{eqnarray*}
In the last few steps, we used the strong Markov property at $\tau
^{\ast}(u+za(u))$, equated $P(X_{tr(u)-v}\in(u+y-\Delta,u+y])$ with
$P(X_{tr(u)-v}^{\ast}\in
(-u-y+\Delta,-u-y])$, and used (\ref{llb}). Since $c(tr(u))\sim
t^{1/\overline{\gamma}}c(r(u))=t^{1/\overline{\gamma}}a(u)$, this
together with \eqref{6.8a} gives (\ref{f}).

Now for case \textup{(i)}, with $a(u)=u$
and $P(C>z)=(1+z)^{-\beta}$,
%
%e6.11 #&#
\begin{eqnarray}
\label{pne} &&P^{(u)} \bigl(\tau_u\in r(u)\,\mathrm{d}t,
Z^{(u)}\in \bigl[za(u),za(u)+\Delta\bigr] \bigr)\nonumber
\\
&&\quad \sim
\frac{ \overline{\Pi}^+_X(u)P(C>z)P(X_{tr(u)}^*\in
[za(u),za(u)+\Delta
])\,\mathrm{d} t} {P(\tau_u<\infty)}\qquad
 \bigl(\mbox{by \eqref{pX} and \eqref{6.6a}}\bigr)\nonumber
\\
&& \quad\sim
\frac{(1+z)^{-\beta}\overline{\Pi
}^+_X(u)h_{t}(z)\Delta
\,\mathrm{d} t } {q^{-1}
\overline{\Pi}_\cH(u)a(u)}\qquad \bigl(\mbox{by \eqref{e1} and \eqref {max}}\bigr)
\\
\nonumber &&\quad \sim
\frac{(1+z)^{-\beta
}qc_{\gamma,\beta}A_{H^*}(u)h_{t}(z)\Delta\,\mathrm{d} t} {a^2(u)}\qquad \bigl(\mbox{by \eqref{c1}, with $A_{H^*}(\infty)=
\infty$}\bigr)
\\
&&\quad \sim\frac{(1+z)^{-\beta
}c(\gamma,\beta)h_{t}(z)\Delta\,\mathrm{d} t}
{\Gamma(1+\gamma)a(u)r(u)}\qquad \bigl(\mbox{by \eqref{g}, and
$c_{\gamma,\beta}=c(\gamma ,\beta)$}\bigr)\nonumber
\\
&& \quad = \frac{(1+z)^{-\beta}\Gamma(\beta)h_{t}(z)\Delta\mathrm
{d} t}{
\Gamma(\beta+\gamma-1)a(u)r(u)}\qquad
\bigl(\mbox{by \eqref{cgbdef}}\bigr).\nonumber
\end{eqnarray}
This gives
\[
\lim_{u\to\infty}a(u)r(u)P^{(u)}\bigl(Z^{(u)}\in
\bigl(za(u),za(u)+\Delta\bigr],\tau _u\in r(u)\,\mathrm{d} t\bigr)=
h_{t}(z)f(z)\Delta\,\mathrm{d} t,
\]
where $f(\cdot)$ is as defined in \eqref{fdef}, and proves \eqref{Z} for
case \textup{(i)}.

In case \textup{(ii)}, we get from \eqref{est}
%\[
%\frac{a(u)}{A_{H^*}(u)} \sim q\Gamma(\gamma+1)r(u)
%\mbox{ and }
%\frac{a(u)\overline{F}(u)}{\overline{%
%\Pi_{\mathcal{H}} }(u)}=\frac{1}{\Gamma(2-\gamma)},
%\]%
%and since $\Gamma(2-\gamma)\Gamma(\gamma+1)k_{\gamma}=1$ we easily
%deduce that
\[
P^{(u)}\bigl(Z^{(u)}\in\bigl(za(u),za(u)+\Delta\bigr],\tau_u
\in r(u)\,\mathrm{d} t\bigr) \sim \mathrm{e}^{-z}\frac{h_{t}(z)\Delta}{r(u)a(u)}\,\mathrm{d} t=
\frac
{h_{t}(z)f(z)\Delta
}{r(u)a(u)}\,\mathrm{d} t,
\]
and (\ref{Z}) is established in this case.

Notice also that, since $h_{t}(\cdot)$ vanishes on the negative half-line,
the previous estimates show that $P^{(u)}(-Z^{(u)}\in
(za(u),za(u)+\Delta
],\tau_u\in r(u)\,\mathrm{d} t)/\mathrm{d} t$ is uniformly $\mathrm{o}((r(u)a(u))^{-1})$
for $z\in[\Delta_0,\Delta_1]$ and $t\in[0,T_0]$.

We have now proved \eqref{Z}. It remains to prove (\ref{theta}).

For $k\geq2$, we assume first that $z_{1}<z_{2}<\cdots<z_{k}$ and
write (\ref{theta}) as
\[
\bigl(a(u)\bigr)^{k}r(u)P^{(u)} \Biggl( \bigcap
_{i=1}^{k}C_{i}\cap B \Biggr) =\theta
_{k}(z_{1},z_{2},\ldots, z_{k},t)
\Biggl( \prod_{i=1}^{k}\Delta
_{i}+\mathrm{o}(1) \Biggr) \,\mathrm{d} t,
\]
where %%% again writing $tr(u)=r$
%e6.12 #&#
\begin{eqnarray}
\label{ee} C_{i}&:= &\bigl\{ X^{\ast}\bigl(s_{i}tr(u)
\bigr)\in\bigl(z_{i}a(u),z_{i}a(u)+\Delta _{i}\bigr]
\bigr\},\qquad i=1,2,\ldots,k\quad \mathrm{and}
\nonumber
\\[-8pt]
\\[-8pt]
\nonumber
 B&:=& \bigl\{ \tau _u\in r(u)
\,\mathrm{d} t \bigr\}.
\end{eqnarray}
As in the lines leading up to \eqref{e1}, we have
%
%e6.13 #&#
\begin{equation}
\label{e} P \Biggl( \bigcap_{i=1}^{k}C_{i}
\cap B \Biggr) \sim P \Biggl( \bigcap_{i=1}^{k}C_{i}
\cap\{\overline{X}_{tr(u)}\leq u\} \Biggr) \overline {
\Pi}_{X}^{+}\bigl(u+z_{k}a(u)\bigr)\,\mathrm{d}t.
\end{equation}
The event in brackets on the RHS coincides with $\bigcap_{i=1}^{k}\widetilde{C}_{i}$ where
\[
\widetilde{C}_{i}:= \Bigl\{ X^{\ast}(s_{i}r)\in
\bigl(z_{i}a(u),z_{i}a(u)+\Delta_{i}\bigr],\sup
_{rs_{i-1}<v\leq
rs_{i}}X_{v}\leq u \Bigr\}
\]
and we set $r:=tr(u)$. Note that each $r(u)(s_{i}-s_{i-1})\rightarrow
\infty
$ uniformly in $ i=1,2,\ldots,k $ as $u\rightarrow\infty$. So by the
Markov property and stationarity we find that $P ( \bigcap_{i=1}^{k}\widetilde{C}_{i} )$ is equal to
\begin{eqnarray*}
&&\int_{a(u)z_{k-1}}^{a(u)z_{k-1}+\Delta_{k-1}} P\Bigl (X^{\ast
}(rs_{k})\in \bigl(z_ka(u),z_ka(u)+\Delta_{k}\bigr], \sup_{rs_{k-1}<v\leq rs_{k}}X_{v}\leq u\mid X^{\ast
}(rs_{k-1})=y
\Bigr)\\
&&\hspace*{42pt}\qquad{} \times P \Biggl( \bigcap_{i=1}^{k-1}
\tilde{C}_{i}, X^{\ast
}(rs_{k-1})\in \mathrm{d}y
\Biggr)
\\
&&\quad=  \int_{a(u)z_{k-1}}^{a(u)z_{k-1}+\Delta_{k-1}} P \bigl(X^{\ast
}\bigl(r(s_{k}-s_{k-1})\bigr)
\in\bigl(z_ka(u)-y,z_ka(u)-y+\Delta_{k}\bigr],\\
&&\hspace*{98pt}\overline{X}_{r(s_{k}-s_{k-1})}\leq u-y \bigr) \\
&&\hspace*{62pt}\qquad{}\times P \Biggl( \bigcap
_{i=1}^{k-1}\tilde{C}_{i},X^{\ast
}(rs_{k-1})\in \mathrm{d}y \Biggr)
\\
&&\quad=  \frac{\Delta_{k}}{a(u)} \bigl( h_{t(s_{k}-s_{k-1)}} \bigl( (z_{k}-z_{k-1})
\bigr) +\mathrm{o}(1) \bigr) \times P^{(u)} \Biggl(\bigcap
_{i=1}^{k-1}\widetilde C_{i} \Biggr),
\end{eqnarray*}
where the last line uses the result for $k=1$ in \eqref{Z}. Repeating
this argument, a further
$k-1$ times gives
\[
P \Biggl( \bigcap_{i=1}^{k}\widetilde
C_i \Biggr) =\bigl(a(u)\bigr)^{-k}\prod
_{i=1}^{k}\Delta_{i} \Biggl( \prod
_{i=1}^{k}h_{t(s_{i}-s_{i-1})}(z_{i}-z_{i-1})+\mathrm{o}(1)
\Biggr),
\]
and the result then follows from \eqref{e} and the previous calculation.
Clearly, if any $z_{i}\leq z_{i-1}$ the calculation is still valid, but
the above product vanishes.
\end{pf}

Using this local result and Lemma~\ref{Lovun}, we easily obtain convergence
of the finite-dimensional distributions, as claimed in part 3.

Now argue as follows. Equation~\eqref{Z} implies that $Z^{(u)}/a(u)$ has a
proper limiting distribution under $P^{(u)}$. By Lemma~\ref{Lovun},
this means that $(Z^{(u)}/a(u),\mathrm{O}^{(u)}/a(u))$ has a proper limiting
distribution under $P^{(u)}$, thus, in
particular, $\mathrm{O}^{(u)}/a(u)$ has a proper limiting distribution under
$P^{(u)}$. From
Proposition~\ref{nasc}, we then deduce Properties 1\textup{(a)} and 1\textup{(b)}, and the
proof of Theorem~\ref{new} is completed by repeating the tightness argument
of the previous section, almost word for word.\qed

%re6.1 #&#
\begin{rem}\label{6.1}
Assumption \eqref{LHsub}, that $\mathcal{H}\in\mathcal{S}
$, is only needed for application of Proposition~\ref{nasc}, where it is
used in effect to deduce that $\overline{\Pi}_\cH(u)\sim qP(\tau
_u<\infty)$
via \eqref{max}. We could replace assumption \eqref{LHsub} with the
assumption $\overline{\Pi}_\cH(u)\sim qP(\tau_u<\infty)$
throughout. But
general necessary and sufficient conditions for the latter in terms of more
basic quantities are currently not known.

Further note that $\overline{\Pi}_\cH(u)$ is not asymptotically
equivalent to the more basic quantity $\overline{\Pi}_X^+(u)$ in our
situation. Vigon's ``\'{e}quation amicale inver\'{s}ee'' is
%
%e6.14 #&#
\begin{equation}
\label{Vig5} \overline{\Pi}_{\mathcal{H}}(u)=\int_{(0,\infty)}
\overline{\Pi }_X^+(y+u)G^*(\mathrm{d} y)
\end{equation}
(recall that $G^{\ast}$ is the renewal measure in the down-going
ladder height process $H^*$, see \eqref{m}). Under the assumption
$\lim_{t\to\infty} X_{t} =
-\infty$ a.s., we have
$G^*(\infty)=\infty$, and it is not hard to show from
\eqref{Vig5} that either $\overline{\Pi}_\cH\in\mathcal{L}$ (see
\eqref{classL}, or $\overline{\Pi}_X^+\in\mathcal{L}$ implies
$\overline{\Pi}_\cH(u)/\overline{\Pi}_X^+(u)\to\infty$.

In general, a sufficient condition for $\overline{\Pi}_\cH\in
\mathcal
{S}$ is $\overline{\Pi}_X^+\in\mathcal{D\cap L}$, where $\mathcal{D}$
is the class of dominatedly varying functions; that is, those for which
$\limsup_{x\to\infty}\overline{\Pi}_X^+(x/2)/\overline{\Pi
}_X^+(x)<\infty$;
see, for example, Foss, Korshunov and Zachary \cite{FKZ}, page~11. So we
can replace Assumption \eqref{LHsub} by $\overline{\Pi}_X^+\in
\mathcal
{D}\cap\mathcal{L}$ throughout.
In particular, $\overline{\Pi}_X^+\in\mathcal{D}$ if $\overline
{\Pi
}_X^+$ is regularly varying with index $-\alpha$ for $\alpha\ge0$.

Further connections between $\overline{\Pi}_\cH$ and $\overline{\Pi}_X$
are in Proposition~5.4 of Kl\"{u}ppelberg, Kyprianou and Maller \cite
{kkm} and the related discussion.
\end{rem}

%s7 #&#
\section{Random walks and compound Poisson processes}\label{s6}

We can specialize our results to the case that $X$ is a compound
Poisson process of the form $X_{t}=S_{N_{t}}$, where $(S_{n},n\geq0)$
is a
random walk and $(N_{t},t\geq0)$ is an independent Poisson counting process
of unit rate. Then, writing $Z_{n}$ and $Z_{n}^{\ast}$ for the $n$th strict
increasing and weak decreasing ladder heights in $S$, we have also that
$H_{t}=Z_{N_{t}}$ and $H_{t}^{\ast}=Z_{N_{t}}^{\ast}$ for all $t\geq0$.
Then our basic assumptions, (\ref{L-inf}) and (\ref{LHsub}) are
equivalent to
\[
S_{n}\overset{\mathrm{a.s.}} {\rightarrow}-\infty\quad\mbox{and}\quad J\in\mathcal{S},
\]
where $J(\mathrm{d} x)=P(Z_{1}\in\mathrm{d} x\mid Z_{1}\in(0,\infty
))$. It is
also clear that, with $\tau^{S}(u):=\inf\{n\dvt S_{n}>u\}$, we have the
identity
\[
\tau_u=\sum_{1}^{\tau^{S}(u)}e_{i},
\]
where the $e_{i}$ are i.i.d. $\operatorname{Exp}(1)$ random variables. Clearly, the
event $\{\tau_u<\infty\}$ coincides a.s. with the event $\{\tau
^{S}(u)<\infty\} $, so $P^{(u)}(\cdot)$ has an unambiguous meaning and, furthermore, it is
straightforward to show that for any $r(u)\rightarrow\infty$ as
$u\rightarrow
\infty$, the statements
\[
r(u)P^{(u)}\bigl(\tau^{S}(u)=\bigl[tr(u)\bigr]\bigr)
\rightarrow g(t)
\]
and
\[
r(u)P^{(u)}\bigl(\tau_u \in r(u)\,\mathrm{d} t\bigr)
\rightarrow g(t)\,\mathrm{d} t
\]
are equivalent. Also the spatial quantities $Z_{S}^{(u)}:=S^{\ast
}(\tau
^{S}(u))$ and $\mathrm{O}_{S}^{(u)}:=S(\tau^{S}(u))-u$ coincide with $Z^{(u)}$
and $\mathrm{O}^{(u)}$.

We claim that this allows us to deduce versions of Theorems~\ref{new}
and~\ref{old} for random walks, with very minor changes. Specifically,
if $F$ is
the distribution of $S_{1}$ and we replace $\Pi$ and $\Pi_{\mathcal
{H}} $
in those results by $F$ and $J$, then Theorem~\ref{old} requires only
replacing $g^{(u)}(tr(u))$ by $P^{(u)} (\tau^{S}(u)=[tr(u)] )$,
and Theorem~\ref{new} requires only an analogous change to (\ref{Z}).

Alternatively, we can prove the random walk results by repeating the L\'
{e}vy process proof, with appropriate changes.
We refer to Borovkov and Borovkov \cite{BandB} for general results on
heavy-tailed
random walks.

%re7.1 #&#
\begin{rem}\label{7.1}
An alternative approach to our proofs, suggested by a referee, based on
``the principle of a single large jump'' (developed in Asmussen and Foss
\cite{AF} for a more general
setting and then considered in Chapter~5, Section~13 of Foss, Korshunov
and Zachary \cite{FKZ} for random walks), may provide a shorter and
more intuitive treatment. However, extending these techniques to the L\'
{e}vy process situation and dealing with the infinite mean case is not
straightforward, and it is not clear that
this approach would deliver the local results or the if and only if
conditions which we establish.
\end{rem}

% zodis Acknowledgments" paliekamas pagal autoriu
\section*{Acknowledgements}
Research partially supported by ARC Grant DP1092502.
%%%%%%%%%%%%%%%%%%%%%%%%%%%%%%%%%%%%%%%%%%%%%%%%%%%%%%%%%%%%%%%%%%%%%%%%%%%%%%%%%%%%%%%%%%%%%%%%%%%%%%%%%%
% Bibliography
%%%%%%%%%%%%%%%%%%%%%%%%%%%%%%%%%%%%%%%%%%%%%%%%%%%%%%%%%%%%%%%%%%%%%%%%%%%%%%%%%%%%%%%%%%%%%%%%%%%%%%%%%%

% imsref loaded by akundreckaite, 2015-04-29 11:06:49

%\begin{appendix}
%\section{}
%\end{appendix}

%\begin{supplement}%[id=suppA]
%\sname{Supplement A}
%\stitle{}
%\slink[doi]{10.3150/00-BEJXXXXSUPP} %[doi,text={...}] - jei reikia
%suskaldyti doi
%\sdatatype{.pdf}
%\sfilename{BEJ000\_supp.pdf}
%\sdescription{}
%\end{supplement}

%\begin{thebibliography}{00}
%\bibitem{r1}
%\bibitem{r1}
%\end{thebibliography}

\printhistory

\begin{thebibliography}{24}
% pybtex-1.31. Style name=bej, version=1.42, label_style=nolabel, sorting_style=complex, cfg=None, language=None.


%b1 ###
%b1 #&#
\bibitem{AF}
\begin{barticle}[mr]
\bauthor{\bsnm{Asmussen},~\bfnm{S{\o}ren}\binits{S.}} \AND
\bauthor{\bsnm{Foss},~\bfnm{Sergey}\binits{S.}}
(\byear{2014}).
\btitle{On exceedance times for some processes with dependent increments}.
\bjournal{J. Appl. Probab.}
\bvolume{51}
\bpages{136--151}.
\bid{doi={10.1239/jap/1395771419}, issn={0021-9002}, mr={3189447}}
\end{barticle}
%

\bptok{imsref}%
% NOT OUTPUTTED:
%   number = 1
%   doi = http://dx.doi.org/10.1239/jap/1395771419
%   fjournal = Journal of Applied Probability
\endbibitem

%b2 ###
%b2 #&#
\bibitem{AK}
\begin{barticle}[mr]
\bauthor{\bsnm{Asmussen},~\bfnm{S{\o}ren}\binits{S.}} \AND
\bauthor{\bsnm{Kl{\"u}ppelberg},~\bfnm{Claudia}\binits{C.}}
(\byear{1996}).
\btitle{Large deviations results for subexponential tails, with applications to insurance risk}.
\bjournal{Stochastic Process. Appl.}
\bvolume{64}
\bpages{103--125}.
\bid{doi={10.1016/S0304-4149(96)00087-7}, issn={0304-4149}, mr={1419495}}
\end{barticle}
%

\bptok{imsref}%
% NOT OUTPUTTED:
%   number = 1
%   doi = http://dx.doi.org/10.1016/S0304-4149(96)00087-7
%   coden = STOPB7
%   fjournal = Stochastic Processes and their Applications
\endbibitem

%b3 ###
%b3 #&#
\bibitem{bert}
\begin{bbook}[mr]
\bauthor{\bsnm{Bertoin},~\bfnm{Jean}\binits{J.}}
(\byear{1996}).
\btitle{L\'evy Processes}.
\bseries{Cambridge Tracts in Mathematics}
\bvolume{121}.
\blocation{Cambridge}:
\bpublisher{Cambridge Univ. Press}.
\bid{mr={1406564}}
\end{bbook}
%

\bptok{imsref}%
% NOT OUTPUTTED:
%   isbn = 0-521-56243-0
%   fpage = x+265
\endbibitem

%b4 ###
%b4 #&#
\bibitem{BD}
\begin{barticle}[mr]
\bauthor{\bsnm{Bertoin},~\bfnm{J.}\binits{J.}} \AND
\bauthor{\bsnm{Doney},~\bfnm{R.~A.}\binits{R.A.}}
(\byear{1994}).
\btitle{Cram\'er's estimate for L\'evy processes}.
\bjournal{Statist. Probab. Lett.}
\bvolume{21}
\bpages{363--365}.
\bid{doi={10.1016/0167-7152(94)00032-8}, issn={0167-7152}, mr={1325211}}
\end{barticle}
%

\bptok{imsref}%
% NOT OUTPUTTED:
%   number = 5
%   doi = http://dx.doi.org/10.1016/0167-7152(94)00032-8
%   coden = SPLTDC
%   fjournal = Statistics \& Probability Letters
\endbibitem

%b5 ###
%b5 #&#
\bibitem{BGT}
\begin{bbook}[mr]
\bauthor{\bsnm{Bingham},~\bfnm{N.~H.}\binits{N.H.}},
\bauthor{\bsnm{Goldie},~\bfnm{C.~M.}\binits{C.M.}} \AND
\bauthor{\bsnm{Teugels},~\bfnm{J.~L.}\binits{J.L.}}
(\byear{1987}).
\btitle{Regular Variation}.
\bseries{Encyclopedia of Mathematics and Its Applications}
\bvolume{27}.
\blocation{Cambridge}:
\bpublisher{Cambridge Univ. Press}.
\bid{doi={10.1017/CBO9780511721434}, mr={0898871}}
\end{bbook}
%

\bptok{imsref}%
% NOT OUTPUTTED:
%   doi = http://dx.doi.org/10.1017/CBO9780511721434
%   isbn = 0-521-30787-2
%   fpage = xx+491
\endbibitem

%b6 ###
%b6 #&#
\bibitem{BandB}
\begin{bbook}[mr]
\bauthor{\bsnm{Borovkov},~\bfnm{A.~A.}\binits{A.A.}} \AND
\bauthor{\bsnm{Borovkov},~\bfnm{K.~A.}\binits{K.A.}}
(\byear{2008}).
\btitle{Asymptotic Analysis of Random Walks: Heavy-Tailed Distributions}.
\bseries{Encyclopedia of Mathematics and Its Applications}
\bvolume{118}.
\blocation{Cambridge}:
\bpublisher{Cambridge Univ. Press}.
%\bnote{, Translated from the Russian by O. B. Borovkova}.
\bid{doi={10.1017/CBO9780511721397}, mr={2424161}}
\end{bbook}
%

\bptok{imsref}%
% NOT OUTPUTTED:
%   doi = http://dx.doi.org/10.1017/CBO9780511721397
%   isbn = 978-0-521-88117-3
%   fpage = xxx+625
\endbibitem

%b7 ###
%b7 #&#
\bibitem{HaanFerr}
\begin{bbook}[mr]
\bauthor{\bparticle{de} \bsnm{Haan},~\bfnm{Laurens}\binits{L.}} \AND
\bauthor{\bsnm{Ferreira},~\bfnm{Ana}\binits{A.}}
(\byear{2006}).
\btitle{Extreme Value Theory: An Introduction}.
\bseries{Springer Series in Operations Research and Financial Engineering}.
\blocation{New York}:
\bpublisher{Springer}.
\bid{doi={10.1007/0-387-34471-3}, mr={2234156}}
\end{bbook}
%

\bptok{imsref}%
% NOT OUTPUTTED:
%   doi = http://dx.doi.org/10.1007/0-387-34471-3
%   isbn = 978-0-387-23946-0; 0-387-23946-4
%   fpage = xviii+417
\endbibitem

%b8 ###
%b8 #&#
\bibitem{DFK}
\begin{barticle}[mr]
\bauthor{\bsnm{Denisov},~\bfnm{Denis}\binits{D.}},
\bauthor{\bsnm{Foss},~\bfnm{Serguei}\binits{S.}} \AND
\bauthor{\bsnm{Korshunov},~\bfnm{Dima}\binits{D.}}
(\byear{2004}).
\btitle{Tail asymptotics for the supremum of a random walk when the mean is not finite}.
\bjournal{Queueing Syst.}
\bvolume{46}
\bpages{15--33}.
\bid{doi={10.1023/B:QUES.0000021140.87161.9c}, issn={0257-0130}, mr={2072274}}
\end{barticle}
%

\bptok{imsref}%
% NOT OUTPUTTED:
%   number = 1-2
%   doi = http://dx.doi.org/10.1023/B:QUES.0000021140.87161.9c
%   fjournal = Queueing Systems. Theory and Applications
\endbibitem

%b9 ###
%b9 #&#
\bibitem{doneystf}
\begin{bbook}[mr]
\bauthor{\bsnm{Doney},~\bfnm{Ronald~A.}\binits{R.A.}}
(\byear{2007}).
\btitle{Fluctuation Theory for L\'evy Processes}.
\bseries{Lecture Notes in Math.}
\bvolume{1897}.
\blocation{Berlin}:
\bpublisher{Springer}.
%\bnote{Lectures from the 35th Summer School on Probability Theory held in Saint-Flour, July 6--23, 2005, Edited and with a foreword by Jean Picard}.
\bid{mr={2320889}}
\end{bbook}
%

\bptok{imsref}%
% NOT OUTPUTTED:
%   isbn = 978-3-540-48510-0; 3-540-48510-4
%   fpage = x+147
\endbibitem

%b10 ###
%b10 #&#
\bibitem{DK}
\begin{barticle}[mr]
\bauthor{\bsnm{Doney},~\bfnm{R.~A.}\binits{R.A.}} \AND
\bauthor{\bsnm{Kyprianou},~\bfnm{A.~E.}\binits{A.E.}}
(\byear{2006}).
\btitle{Overshoots and undershoots of L\'evy processes}.
\bjournal{Ann. Appl. Probab.}
\bvolume{16}
\bpages{91--106}.
\bid{doi={10.1214/105051605000000647}, issn={1050-5164}, mr={2209337}}
\end{barticle}
%

\bptok{imsref}%
% NOT OUTPUTTED:
%   number = 1
%   doi = http://dx.doi.org/10.1214/105051605000000647
%   fjournal = The Annals of Applied Probability
\endbibitem

%b11 ###
%b11 #&#
\bibitem{DM2002}
\begin{barticle}[mr]
\bauthor{\bsnm{Doney},~\bfnm{R.~A.}\binits{R.A.}} \AND
\bauthor{\bsnm{Maller},~\bfnm{R.~A.}\binits{R.A.}}
(\byear{2002}).
\btitle{Stability and attraction to normality for L\'evy processes at zero and at infinity}.
\bjournal{J. Theoret. Probab.}
\bvolume{15}
\bpages{751--792}.
\bid{doi={10.1023/A:1016228101053}, issn={0894-9840}, mr={1922446}}
\end{barticle}
%

\bptok{imsref}%
% NOT OUTPUTTED:
%   number = 3
%   doi = http://dx.doi.org/10.1023/A:1016228101053
%   coden = JTPREO
%   fjournal = Journal of Theoretical Probability
\endbibitem

%b12 ###
%b12 #&#
\bibitem{DR2012}
\begin{barticle}[mr]
\bauthor{\bsnm{Doney},~\bfnm{R.~A.}\binits{R.A.}} \AND
\bauthor{\bsnm{Rivero},~\bfnm{V.}\binits{V.}}
(\byear{2013}).
\btitle{Asymptotic behaviour of first passage time distributions for L\'evy processes}.
\bjournal{Probab. Theory Related Fields}
\bvolume{157}
\bpages{1--45}.
\bid{doi={10.1007/s00440-012-0448-x}, issn={0178-8051}, mr={3101839}}
\bptnote{check volume}%
\end{barticle}
%

\bptok{imsref}%
% NOT OUTPUTTED:
%   number = 1-2
%   doi = http://dx.doi.org/10.1007/s00440-012-0448-x
%   fjournal = Probability Theory and Related Fields
\endbibitem

%b13 ###
%b13 #&#
\bibitem{EGV}
\begin{barticle}[mr]
\bauthor{\bsnm{Embrechts},~\bfnm{Paul}\binits{P.}},
\bauthor{\bsnm{Goldie},~\bfnm{Charles~M.}\binits{C.M.}} \AND
\bauthor{\bsnm{Veraverbeke},~\bfnm{No{\"e}l}\binits{N.}}
(\byear{1979}).
\btitle{Subexponentiality and infinite divisibility}.
\bjournal{Z. Wahrsch. Verw. Gebiete}
\bvolume{49}
\bpages{335--347}.
\bid{doi={10.1007/BF00535504}, issn={0044-3719}, mr={0547833}}
\end{barticle}
%

\bptok{imsref}%
% NOT OUTPUTTED:
%   number = 3
%   doi = http://dx.doi.org/10.1007/BF00535504
%   fjournal = Zeitschrift f\"ur Wahrscheinlichkeitstheorie und Verwandte Gebiete
\endbibitem

%b14 ###
%b14 #&#
\bibitem{EKM}
\begin{bbook}[mr]
\bauthor{\bsnm{Embrechts},~\bfnm{Paul}\binits{P.}},
\bauthor{\bsnm{Kl{\"u}ppelberg},~\bfnm{Claudia}\binits{C.}} \AND
\bauthor{\bsnm{Mikosch},~\bfnm{Thomas}\binits{T.}}
(\byear{1997}).
\btitle{Modelling Extremal Events: For Insurance and Finance}.
\bseries{Applications of Mathematics (New York)}
\bvolume{33}.
\blocation{Berlin}:
\bpublisher{Springer}.
\bid{doi={10.1007/978-3-642-33483-2}, mr={1458613}}
\end{bbook}
%

\bptok{imsref}%
% NOT OUTPUTTED:
%   doi = http://dx.doi.org/10.1007/978-3-642-33483-2
%   isbn = 3-540-60931-8
%   fpage = xvi+645
\endbibitem

%b15 ###
%b15 #&#
\bibitem{FKZ}
\begin{bbook}[mr]
\bauthor{\bsnm{Foss},~\bfnm{Sergey}\binits{S.}},
\bauthor{\bsnm{Korshunov},~\bfnm{Dmitry}\binits{D.}} \AND
\bauthor{\bsnm{Zachary},~\bfnm{Stan}\binits{S.}}
(\byear{2013}).
\btitle{An Introduction to Heavy-Tailed and Subexponential Distributions},
\bedition{2nd} ed.
\bseries{Springer Series in Operations Research and Financial Engineering}.
\blocation{New York}:
\bpublisher{Springer}.
\bid{doi={10.1007/978-1-4614-7101-1}, mr={3097424}}
\end{bbook}
%

\bptok{imsref}%
% NOT OUTPUTTED:
%   doi = http://dx.doi.org/10.1007/978-1-4614-7101-1
%   isbn = 978-1-4614-7100-4; 978-1-4614-7101-1
%   fpage = xii+157
\endbibitem

%b16 ###
%b16 #&#
\bibitem{GR}
\begin{barticle}[mr]
\bauthor{\bsnm{Goldie},~\bfnm{Charles~M.}\binits{C.M.}} \AND
\bauthor{\bsnm{Resnick},~\bfnm{Sidney}\binits{S.}}
(\byear{1988}).
\btitle{Distributions that are both subexponential and in the domain of attraction of an extreme-value distribution}.
\bjournal{Adv. in Appl. Probab.}
\bvolume{20}
\bpages{706--718}.
\bid{doi={10.2307/1427356}, issn={0001-8678}, mr={0967994}}
\end{barticle}
%

\bptok{imsref}%
% NOT OUTPUTTED:
%   number = 4
%   doi = http://dx.doi.org/10.2307/1427356
%   coden = AAPBBD
%   fjournal = Advances in Applied Probability
\endbibitem

%%b17 ###
%%b17 #&#
%\bibitem{GM4}
%\begin{barticle}[mr]
%\bauthor{\bsnm{Griffin},~\bfnm{Philip~S.}\binits{P.S.}} \AND
%\bauthor{\bsnm{Maller},~\bfnm{Ross~A.}\binits{R.A.}}
%(\byear{2011}).
%\btitle{The time at which a L\'evy process creeps}.
%\bjournal{Electron. J. Probab.}
%\bvolume{16}
%\bpages{2182--2202}.
%\bid{doi={10.1214/EJP.v16-945}, issn={1083-6489}, mr={2861671}}
%\bptnote{check pages}%
%\end{barticle}
%%
%
%\bptok{imsref}%
%% NOT OUTPUTTED:
%%   doi = http://dx.doi.org/10.1214/EJP.v16-945
%%   fjournal = Electronic Journal of Probability
%\endbibitem

%b18 ###
%b18 #&#
\bibitem{kk}
\begin{barticle}[mr]
\bauthor{\bsnm{Kl{\"u}ppelberg},~\bfnm{C.}\binits{C.}} \AND
\bauthor{\bsnm{Kyprianou},~\bfnm{A.~E.}\binits{A.E.}}
(\byear{2006}).
\btitle{On extreme ruinous behaviour of L\'evy insurance risk processes}.
\bjournal{J. Appl. Probab.}
\bvolume{43}
\bpages{594--598}.
\bid{doi={10.1239/jap/1152413744}, issn={0021-9002}, mr={2248586}}
\end{barticle}
%

\bptok{imsref}%
% NOT OUTPUTTED:
%   number = 2
%   doi = http://dx.doi.org/10.1239/jap/1152413744
%   coden = JPRBAM
%   fjournal = Journal of Applied Probability
\endbibitem

%b19 ###
%b19 #&#
\bibitem{kkm}
\begin{barticle}[mr]
\bauthor{\bsnm{Kl{\"u}ppelberg},~\bfnm{Claudia}\binits{C.}},
\bauthor{\bsnm{Kyprianou},~\bfnm{Andreas~E.}\binits{A.E.}} \AND
\bauthor{\bsnm{Maller},~\bfnm{Ross~A.}\binits{R.A.}}
(\byear{2004}).
\btitle{Ruin probabilities and overshoots for general L\'evy insurance risk processes}.
\bjournal{Ann. Appl. Probab.}
\bvolume{14}
\bpages{1766--1801}.
\bid{doi={10.1214/105051604000000927}, issn={1050-5164}, mr={2099651}}
\end{barticle}
%

\bptok{imsref}%
% NOT OUTPUTTED:
%   number = 4
%   doi = http://dx.doi.org/10.1214/105051604000000927
%   fjournal = The Annals of Applied Probability
\endbibitem

%b20 ###
%b20 #&#
\bibitem{Pakes}
\begin{barticle}[mr]
\bauthor{\bsnm{Pakes},~\bfnm{Anthony~G.}\binits{A.G.}}
(\byear{2004}).
\btitle{Convolution equivalence and infinite divisibility}.
\bjournal{J. Appl. Probab.}
\bvolume{41}
\bpages{407--424}.
\bid{issn={0021-9002}, mr={2052581}}
\end{barticle}
%

\bptok{imsref}%
% NOT OUTPUTTED:
%   number = 2
%   coden = JPRBAM
%   fjournal = Journal of Applied Probability
\endbibitem

%b21 ###
%b21 #&#
\bibitem{Pakes2}
\begin{barticle}[mr]
\bauthor{\bsnm{Pakes},~\bfnm{Anthony~G.}\binits{A.G.}}
(\byear{2007}).
\btitle{Convolution equivalence and infinite divisibility: Corrections and corollaries}.
\bjournal{J.~Appl. Probab.}
\bvolume{44}
\bpages{295--305}.
\bid{doi={10.1239/jap/1183667402}, issn={0021-9002}, mr={2340199}}
\end{barticle}
%

\bptok{imsref}%
% NOT OUTPUTTED:
%   number = 2
%   doi = http://dx.doi.org/10.1239/jap/1183667402
%   coden = JPRBAM
%   fjournal = Journal of Applied Probability
\endbibitem

%b22 ###
%b22 #&#
\bibitem{resnick}
\begin{bbook}[mr]
\bauthor{\bsnm{Resnick},~\bfnm{Sidney~I.}\binits{S.I.}}
(\byear{1987}).
\btitle{Extreme Values, Regular Variation, and Point Processes}.
\bseries{Applied Probability. A~Series of the Applied Probability Trust}
\bvolume{4}.
\blocation{New York}:
\bpublisher{Springer}.
\bid{doi={10.1007/978-0-387-75953-1}, mr={0900810}}
\end{bbook}
%

\bptok{imsref}%
% NOT OUTPUTTED:
%   doi = http://dx.doi.org/10.1007/978-0-387-75953-1
%   isbn = 0-387-96481-9
%   fpage = xii+320
\endbibitem

%b23 ###
%b23 #&#
\bibitem{sato}
\begin{bbook}[mr]
\bauthor{\bsnm{Sato},~\bfnm{Ken-iti}\binits{K.-i.}}
(\byear{1999}).
\btitle{L\'evy Processes and Infinitely Divisible Distributions}.
\bseries{Cambridge Studies in Advanced Mathematics}
\bvolume{68}.
\blocation{Cambridge}:
\bpublisher{Cambridge Univ. Press}.
%\bnote{Translated from the 1990 Japanese original, Revised by the author}.
\bid{mr={1739520}}
\end{bbook}
%

\bptok{imsref}%
% NOT OUTPUTTED:
%   isbn = 0-521-55302-4
%   fpage = xii+486
\endbibitem

%b24 ###
%b24 #&#
\bibitem{Vig}
\begin{barticle}[mr]
\bauthor{\bsnm{Vigon},~\bfnm{Vincent}\binits{V.}}
(\byear{2002}).
\btitle{Votre L\'evy rampe-t-il?}
\bjournal{J. Lond. Math. Soc. (2)}
\bvolume{65}
\bpages{243--256}.
\bid{doi={10.1112/S0024610701002885}, issn={0024-6107}, mr={1875147}}
\end{barticle}
%
\bptok{imsref}%
% NOT OUTPUTTED:
%   number = 1
%   doi = http://dx.doi.org/10.1112/S0024610701002885
%   coden = JLMSAK
%   fjournal = Journal of the London Mathematical Society. Second Series
\endbibitem
\end{thebibliography}
\end{document}